\documentclass[11pt]{article}
\usepackage{amsmath}   
\usepackage{amsfonts}   
\usepackage[all]{xy}

\newcommand{\dlim}{\mathop{\underrightarrow{\rm lim}}\nolimits} 
\newtheorem{ttt}{Theorem}[section] 
\newtheorem{lll}[ttt]{Lemma}
\newtheorem{ccc}[ttt]{Corollary} 
\newtheorem{ppp}[ttt]{Proposition}
 
\newtheorem{itdef}[ttt]{Definition}
\newtheorem{itexa}[ttt]{Example}
\newtheorem{itrem}[ttt]{Remark}
\newenvironment{ddd}{\begin{itdef}\rm}{\end{itdef}}

\newenvironment{rrr}{\begin{itrem}\rm}{\end{itrem}}

\newenvironment{pf}[1][Proof]{
\par\noindent{\em #1}. }{\hfill\framebox(6,6)\par\medskip}

\newcommand{\bbF}{{\mathbb F}}

\newcommand{\bbQ}{{\mathbb Q}}
\newcommand{\bbN}{{\mathbb N}}

\newcommand{\End}{\mathop{\rm End}\nolimits}

\newcommand{\Soc}{\mathop{\rm Soc}\nolimits}
\newcommand{\diag}{\mathop{\rm diag}\nolimits}

\newcommand{\chf}{\mathop{\rm char}\nolimits\bbF}

\newcommand{\ev}{\End V}

\newcommand{\ep}{\varepsilon}
\newcommand{\de}{\delta}
\newcommand{\da}{{\downarrow}}
\newcommand{\Ze}{{\mathbb Z}}
\newcommand{\Se}{{\mathcal S}}
\newcommand{\Ce}{{\mathcal C}}
\newcommand{\Te}{{\mathcal T}}

\begin{document}

\title{
%\vspace{10mm}
Classification of the direct limits of involution simple associative algebras
and the corresponding dimension groups
}

\author{A.\,A.\,Baranov
\\
 Department of Mathematics\\
 University of Leicester\\
 Leicester LE1 7RH\\
{\normalsize {\it e-mail}: ab155@le.ac.uk}
}

\date{}

\maketitle

\begin{abstract} 
A classification of the (countable) direct limits 
of finite dimensional involution simple associative algebras
over an algebraically closed field of arbitrary characteristic
is obtained. This also classifies the corresponding dimension groups. 
The set of invariants consists of two supernatural
numbers and two real parameters.
\end{abstract}

\section{Introduction}
The ground field $\bbF$ is algebraically closed of arbitrary characteristic.
Let $A$ be an associative algebra over $\bbF$ (not necessarily
containing an identity element).
Assume $A$ has an involution, 
that is,
a linear transformation $*$ of $A$ such that
$(a^*)^*=a$ and $(ab)^*=b^*a^*$ for all $a,b\in A$. 
We will sometimes denote this algebra by $(A,*)$ to reflect the fact
that $A$ is an algebra with involution.
Note that our involution is $\bbF$-linear, i.e.
we consider involutions of the first kind only.
The algebra $A$ is called {\em involution simple} if $A^2\ne 0$ and 
it has no non-trivial $*$-invariant ideals. 

We say that an infinite dimensional algebra $A$ 
is {\em locally (semi)simple} if any finite subset of $A$ is contained 
in a finite dimensional (semi)simple subalgebra. 
Note that we do not require $A$ to have an identity element.
If $A$ has an 
involution and these subalgebras can be chosen
involution simple with respect to the inherited involution
then $A$ is called {\em locally involution simple}. 
Observe that $A$ itself is involution simple in that case. 
The aim of this paper is to classify locally involution simple 
associative algebras over $\bbF$ of countable dimension.

Let $A$ be a locally simple 
associative algebra of countable dimension 
over $\bbF$. It follows from the definition that there is a chain 
of simple subalgebras $A_1\subset A_2\subset A_3\subset\dots$ of 
$A$ such that $A=\cup_{i=1}^{\infty} A_i$. One can also view 
$A$ as the direct limit $\dlim A_i$ for the sequence 
\begin{equation}\label{e0}
A_1\to A_2\to A_3\to\dots
\end{equation}
of injective homomorphisms of finite dimensional simple associative algebras $A_i$. 
Since $\bbF$ is algebraically closed, each $A_i$ can be identified
with the algebra $M_{n_i}(\bbF)$ of all $n_i\times n_i$ 
matrices over $\bbF$ for some $n_i$. 
Moreover, each embedding $A_i\to A_{i+1}$ can be written in the 
following matrix form
\begin{equation}\label{e01}
M\mapsto \diag(M,\dots,M,0,\dots,0), \quad M\in M_{n_i}(\bbF).
\end{equation}
Therefore in order to describe locally simple associative algebras
of countable dimension one needs to classify the 
direct limits of the sequences of matrix algebras (\ref{e0}).  
Elliot \cite{elliot} did this in terms of systems of idempotents. 
It has been shown later
that Elliot's invariant can be interpreted in terms of the
$K_0$-functor (see Theorem \ref{dimgr}). As a particular case of our main results we get 
another parametrization 
of these algebras. % (see Theorem~\ref{main}).

Assume now that the algebra $A$ is locally involution simple, i.e.
we have a sequence (\ref{e0}) of involution simple finite dimensional
algebras $A_i$ and $A=\dlim A_i$. Note that all homomorphisms
in  (\ref{e0}) respect the involution
but do not necessarily preserve the identity element. 
It is well known 
that every involution simple finite dimensional
$\bbF$-algebra is either a full matrix algebra
or the direct sum of two isomorphic
matrix algebras. Therefore the combinatorial picture is much more complicated
than in (\ref{e01}). However it is still possible
to provide an explicit parametrization (see our main Theorems~\ref{main}
and \ref{t51}).

In Section \ref{dirlimits} we prove that two locally involution 
simple algebras of the same type (orthogonal, symplectic or special) are isomorphic
if and only if they are isomorphic as associative algebras (Theorem \ref{aareduced}).
This partially reduces the classification problem to locally semisimple
associative algebras. These algebras are normally classified by 
ordered dimension groups (see Theorem \ref{dimgr}). 
However, as it is pointed out in \cite{br}, although dimension groups are relatively
easy objects, their isomorphism classes are not, and general classification
is not available. 
Some 
examples of known isomorphism classes of the dimension groups can be found
in \cite{br}. 
Our main Theorem \ref{main} 
gives complete classification of the dimension groups which correspond 
to the locally involution simple algebras. 
These are the direct limits of the sequences
$\mathbb{Z}^3\to \mathbb{Z}^3\to\dots \to \mathbb{Z}^3\to \dots $
where the embeddings are given by the Bratteli diagrams 
(\ref{bratA2}). 

Another approach ($K$-theoretical in nature) to the classification 
of locally involution simple associative algebras 
can be found within the general theory of compact
group actions on locally semisimple algebras,
see \cite{han,beek}. 
In the case of order $2$ automorphisms this was done by Fack and Mar\'echal
\cite{fm1} (unital embeddings given by the Bratteli diagrams 
(\ref{bratA1}))
and Elliott and Su \cite{el} (in terms of $K$-theoretical invariants).

Our approach uses some technique developed
by Baranov and Zhilinskii for the classification of the diagonal direct
limits of finite dimensional simple Lie algebras over 
an algebraically closed field of characteristic zero \cite{bzh}.
It is shown in \cite{bbz} that there is a natural bijective
correspondence between such Lie algebras and locally involution simple
associative algebras, so the classification should be similar.
Unfortunately the proofs in \cite{bbz,bzh} are very
dependent on characteristic zero and fail to work in positive characteristic. 
In the present paper we provide new, characteristic free, proofs.
However, the case of characteristic 2 still requires special attention
and the classification is slightly different in that case.

Note that
our results do not exhaust the problem of classification
of all involution simple locally finite dimensional associative algebras 
(of countable dimension),
since there are examples of such algebras which are not locally semisimple 
(see \cite{FS,z5}).

\section{Preliminaries}
\label{prelim}

Recall that an associative algebra $A$ with involution is called 
{\em involution simple} if $A^2\ne 0$ and it has no non-trivial $*$-invariant ideals. 
The following is well-known.

\begin{ppp} \label{issum}
Let $A$ be an involution simple associative algebra. 
Then either $A$ is simple as an algebra or $A$ has exactly two non-zero proper 
ideals $B_1$ and $B_2$. Moreover both $B_1$ and $B_2$
are simple algebras, $B_1^*=B_2$ and $A=B_1\oplus B_2$.
\end{ppp}

\begin{pf}
Assume $A$ is not simple. Let $B_1$ be a non-zero proper 
ideal of $A$. Then $B_2=B_1^*$ is also an ideal of $A$. 
Since $B_1+B_2$ and $B_1\cap B_2$ are $*$-invariant ideals
of $A$ and $A$ is involution simple, one has 
$B_1+B_2=A$ and $B_1\cap B_2=0$, i.e. $A=B_1\oplus B_2$. 
Now, if $B$ is a non-zero proper 
ideal of $B_1$ then $B\oplus B^*$ is a non-zero proper 
 $*$-invariant ideal of $B_1\oplus B_2=A$. Therefore 
 $B=B_1$ and both $B_1$ and $B_2$ are simple algebras.
% In particular, $B_i^2\ne0$.

Assume now that $C$ is another non-zero proper 
ideal of $A$. Then by the above argument, $A=C\oplus C^*$. 
If $B_1\subseteq C$ or $B_2\subseteq C$  then it is easy to see that $C=B_1$ or $B_2$.
Assume this is not the case. 
Let $B=B_1\cap C$. Then $B+B^*$ is  a proper $*$-invariant ideal of $A$,
so $B=0$. In particular, $B_1C\subseteq B_1\cap C=0$.
Similarly, $B_2C=0$ and $B_1C^*=B_2C^*=0$. This implies 
$AA=(B_1+B_2)(C+C^*)=0$, which is a contradiction.
\end{pf}

Let $A$ be a finite dimensional associative algebra over $\bbF$ with involution $*$. Assume that $A$ is involution simple. Then by Proposition \ref{issum},
$A$ is either simple or  $A=B\oplus B^*$ the sum of two (anti)isomorphic simple subalgebras. 
Thus, we can identify $A$ with either
$\ev$ or $\ev_1\oplus\ev_2$ for some finite dimensional vector spaces $V$,
$V_1$, and $V_2$ over $\bbF$ with $\dim V_1=\dim V_2$.
By fixing bases of $V$, $V_1$, and $V_2$, one can represent 
the algebras $\ev$ and $\ev_1\oplus\ev_2$ in the matrix forms
$M_n(\bbF)$ and $M_m(\bbF)\oplus M_m(\bbF)$, respectively,
where $n=\dim V$ and $m=\dim V_1=\dim V_2$. 
We say that these are their {\em matrix realizations}.
We say that a matrix realization of $(\ev,*)$ is {\em canonical}
if the involution in the chosen basis has one of the following two forms:
\begin{eqnarray}
X\mapsto X^t, & & \quad X\in M_n(\bbF) \quad 
\mbox{(transpose)}; \label{orth} \\
X\mapsto X^{\tau}, & & \quad X\in M_n(\bbF) \quad 
\mbox{(symplectic transpose).} \label{sympl}
\end{eqnarray}
In the latter case 
$n$ is even and $X^{\tau}=-JX^tJ$ where 
$
J=\diag\left(
\left(
\begin{array}{cc}
0  & 1 \\
-1 & 0
\end{array}
\right),\dots,
\left(
\begin{array}{cc}
0  & 1 \\
-1 & 0
\end{array}
\right)
\right)
$
($n/2$ blocks).
We say that a matrix realization of $(\ev_1 \oplus \ev_2,*)$ is {\em canonical}
if the involution in the chosen basis has the following form:
\begin{equation}\label{atype}
(X_1,X_2)\mapsto (X_2^t,X_1^t), \quad X_1,X_2\in M_m(\bbF).
\end{equation}

It is well known that any finite dimensional
involution simple algebra over an algebraically closed field
has a canonical matrix realization. Indeed, let us first consider
the algebra $\ev$.
Let $b:V\times V\to \bbF$ be a nondegenerate symmetric or
skew-symmetric bilinear form on $V$. For each 
$x\in\ev$ define $\alpha_b(x)$ by the following property
$$
b(\alpha_b(x)v,w)=b(v,xw)\quad{\rm for\ all\ } v,w\in V.
$$
Then the map
$$
\alpha_b:\ev\to\ev
$$
is an involution of the algebra $\ev$, called the {\em adjoint
involution} with respect to $b$.
More exactly we have the following fact.

\begin{ttt}[{\cite[Ch.1, Introduction]{book}}]\label{form}
The map $b\mapsto \alpha_b$ induces a one-to-one correspondence
between the equivalence classes of nondegenerate symmetric and 
skew-symmetric bilinear forms on $V$ modulo 
multiplication by a factor in $\bbF^{\times}$ and involutions
(of the first kind) on $\ev$.
\end{ttt}

Recall that a bilinear form is called {\em alternating} if
$b(v,v)=0$ for all $v\in V$.
Obviously, if $\chf\ne2$, then the form $b$ is alternating if and only if
it is skew-symmetric. If $\chf=2$, then $b$ is alternating if and only if
it is symmetric and for any choice of basis of $V$,
all diagonal entries of the matrix of $b$ are zeros.
An involution $\alpha$ of $\ev$ is called 
{\em symplectic} (resp. {\em orthogonal}) if it is adjoint to
an alternating (resp. symmetric non-alternating) bilinear form on $V$.
Recall that each finite dimensional orthogonal (resp. symplectic)
vector space over an algebraically closed field has
an orthonormal (resp. hyperbolic) basis. That is, the matrix of $b$ in this
basis is either the identity (in the orthogonal case)
or $J$ (see above)
in the symplectic case (see for example 
\cite[Theorems 11.10 and 11.14]{Roman}).
It is easy to see that the adjoint involution in this basis
is canonical, i.e. of the forms (\ref{orth}) and (\ref{sympl}), respectively.
Thus, we get the following well-known fact.

\begin{ppp}\label{invol1}
Let $V$ be a vector space of dimension $n$ over $\bbF$
and let $*$ be an involution of $\ev$. Then the algebra
$(\ev,*)$ has a canonical matrix realization.
\end{ppp}

To prove a similar result for the algebra $\ev_1\oplus\ev_2$, we need
the following simple fact.

\begin{ppp}\label{matinv}
Each involution of the matrix algebra $M_n(\bbF)$ is
of the following form: $X\mapsto CX^tC^{-1}$ where $C$ is an invertible
matrix.
\end{ppp}

\begin{pf}
The matrix transpose $X\mapsto X^t$ is a natural involution
of $M_n(\bbF)$. Thus the map $X\mapsto (X^*)^t$ is an automorphism
of $M_n(\bbF)$. By Skolem-Noether theorem each automorphism
of $M_n(\bbF)$ is inner, i.e. there exists an invertible matrix
$K$ such that $(X^*)^t=K^{-1}XK$. Therefore
$X^*=K^tX^t(K^{-1})^t=K^tX^t(K^t)^{-1}$, as required.
\end{pf}

\begin{ppp}\label{invol2} Let $V_1$ and $V_2$ be vector spaces
of dimension $m$ and let $*$ be an involution of 
the algebra $\ev_1\oplus\ev_2$ such that $(\ev_1)^*=\ev_2$. 
Then for every matrix realization of $\ev_1$ there is
a matrix realization of $\ev_2$ such that the corresponding
matrix realization of $(\ev_1\oplus\ev_2,*)$ is canonical.
\end{ppp}

\begin{pf} 
Fix any matrix realizations of $\ev_1$ and $\ev_2$,
i.e. identify these algebras with the algebra  $M_m(\bbF)$.
Then the map $*:\ev_1\to \ev_2$ gives an involution 
$X\mapsto X^*$ of $M_m(\bbF)$.
By Proposition~\ref{matinv}, $X^*=CX^tC^{-1}$. It remains
to change basis of $V_2$ (i.e. matrix realization
of $\ev_2$), to eliminate $C$.
\end{pf}

Let $A$ be an  involution simple finite dimensional algebra over $\bbF$.
We say that $A$ is of type $\bf S$, or of {\em symplectic} type,
if $A$ is simple
as an algebra and the involution is symplectic. Similarly we define
the {\em orthogonal} type $\bf O$. 
If $A$ is not simple, then we say
that $A$ is of type $\bf A$, or of {\em special} type. 
Note that algebras of type $\bf S$ are not isomorphic
to those of type $\bf O$ (as algebras with involution).
Thus the canonical matrix realizations (as in Propositions~\ref{invol1}
and~\ref{invol2}) 
give a complete classification of finite dimensional involution
simple algebras over an algebraically closed field.

\begin{rrr}\label{other}
We will also use other canonical forms for involutions.
Let $n$ be even. Define the following $n\times n$ matrices:
$$
J_+=\diag\left(
\left(
\begin{array}{cc}
0  & 1 \\
1 & 0
\end{array}
\right),\dots,
\left(
\begin{array}{cc}
0  & 1 \\
1 & 0
\end{array}
\right)
\right)
,\quad
Q_\pm=\left(
\begin{array}{cc}
0  & I \\
\pm I & 0
\end{array}
\right)
$$
where $I$ is the identity $n/2\times n/2$ matrix.
Then $J_+$ and $Q_+$ define nondegenerate bilinear symmetric forms
on the natural $M_n(\bbF)$-module. If $\chf\ne2$, these forms
are non-alternating, so induce orthogonal involutions
on $M_n(\bbF)$: $\tau_+:X\mapsto J_+X^tJ_+$ and $\theta_+:X\mapsto Q_+X^tQ_+$.
In view of Proposition~\ref{invol1}, 
by choosing an appropriate basis, these
involutions can be represented as matrix transpose.
Thus each orthogonal involution (for $\chf\ne2$ and algebras of even degree)
can be represented as $\tau_+$
(resp. $\theta_+$) in a suitable basis. 
Similarly, the involution $\theta_-:X\mapsto -Q_-X^tQ_-$
is symplectic and
each symplectic involution (for any characteristic)
can be represented in this form.
\end{rrr}

The following simple fact will be used later.

\begin{lll}\label{fg}
Let $B$ be a finite dimensional involution simple algebra
of even dimension. Let $e$ be the identity element of $B$.
If $\chf=2$, assume that $B$ is not of type $\bf O$.
Then $B$ has idempotents 
$f$ and $g$ such that $e=f+g$, $fg=gf=0$, and $f^*=g$.
\end{lll}

\begin{pf}
This is obvious if $B$ is of type $\bf A$. Assume that
$B$ is symplectic. Represent 
$B$ as in Proposition~\ref{invol1}.
Then one can 
easily check that $f=\diag(1,0,1,0,\dots,1,0)$
and $g=\diag(0,1,0,1,\dots,0,1)$ are the required idempotents.
If the involution $*$ of $B$ is orthogonal, then by Remark~\ref{other}, 
it can be represented as $X^*=J_+X^tJ_+$, $X\in M_n(\bbF)$. 
Then it is easy to check that the same idempotents $f$ and $g$ 
as in the symplectic case satisfy the required conditions.
\end{pf}

Now we are going to study embeddings of involution simple
algebras, i.e. injective homomorphisms 
$\ep: A_1\to A_2$ 
which respect involution.
We do not require these embeddings to preserve the identity element.
Since the embeddings respect involution we often use the same
symbol ``$*$'' to denote the involution of $A_1$ and $A_2$.
We usually identify $A_1$ with its image $\ep(A_1)$
in $A_2$.
If $A_i$ is of type $\bf A$, we denote by
$B_i$ and $C_i$ its simple components (so $A_i=B_i\oplus C_i$ and 
$B_i\cong C_i$). 
It is convenient to assume $B_i=A_i$ if $A_i$ is of type $\bf S$
or $\bf O$. 
We denote by $e_i$, $f_i$, and $g_i$ the identities
of $A_i$, $B_i$, and $C_i$, respectively. 
Thus $e_i=f_i+g_i$ if $A_i$ is of type $\bf A$,
and $e_i=f_i$ otherwise. Note that $f_i^*=g_i$
if $A_i$ is of type $\bf A$.

Recall that $B_i\cong M_{n_i}(\bbF)$
for some $n_i\in\bbN$. We say that $n_i$ is the {\em degree} of
$A_i$. Denote by $V_i$ the natural $B_i$-module
of dimension $n_i$ and by $W_i$ the natural module for
$C_i$ (if $C_i\ne0$). We consider these modules 
as $A_i$-modules in a natural way. 
If $A_i$ is not of type $\bf A$, we denote by $b_i$ a nondegenerate
bilinear form on $V_i$ corresponding to the involution $*$ on $A_i$
(see Theorem~\ref{form}).

Denote by $T_i$ the trivial
one-dimensional $A_i$-module (with zero action).
Now the restriction 
of the $A_2$-module $V_2$ to $A_1$ is completely
reducible, so can be described as follows.
\begin{equation}\label{sign}
V_2\da A_1=
\underbrace{V_1  \oplus\dots\oplus V_1}_{l}\oplus
\underbrace{W_1\oplus\dots\oplus W_1}_{r}\oplus
\underbrace{T_1  \oplus\dots\oplus T_1}_{z}
\end{equation}
where $l,r,z\in\bbN\cup\{0\}$
and $r=0$ if $A_1$ is not of type $\bf A$.

\begin{ddd}\label{d1}
The triple $(l,r,z)$ in (\ref{sign}) is called the {\em signature} of 
the embedding $\ep:A_1\to A_2$.
\end{ddd}

\begin{rrr}\label{lr}
 If both $A_1$ and $A_2$ are of type $\bf A$, then
the signature depends on the choice 
of the simple components of $A_1$ and $A_2$, e.g. by
swapping $B_1$ and $C_1$ (or $B_2$ and $C_2$, see (\ref{e7}) below), 
the signature 
$(l,r,z)$ is replaced by $(r,l,z)$. Thus we can and will assume that
$l\ge r$.
\end{rrr}

\begin{ddd}\label{canembed}
We say that a homomorphism $\ep: M_{n_1} \to M_{n_2}$ of signature $(l,0,z)$ of two
matrix algebras is {\em canonical} if 
\begin{equation}\label{e82}
\ep(M)=
\diag(\underbrace{M,\dots,M}_{l},
\underbrace{0,\dots,0}_{z}), \quad M\in M_{n_1}(\bbF).
\end{equation}
We say that a homomorphism $\ep: M_{n_1}\oplus M_{n_1} \to M_{n_2}\oplus M_{n_2}$
of signature $(l,r,z)$ is 
{\em canonical} if 
\begin{equation}\label{e7}
\ep(M,N)=
(\diag(\underbrace{M,\dots,M}_{l},
\underbrace{N,\dots,N}_{r},\underbrace{0,\dots,0}_{z}),
\diag(\underbrace{N,\dots,N}_{l},
\underbrace{M,\dots,M}_{r},
\underbrace{0,\dots,0}_{z}))
\end{equation} 
for all $M,N\in M_{n_1}(\bbF)$.

We say that an embedding
$\ep: A_1\to A_2$ 
of finite dimensional involution simple algebras over $\bbF$ of the same type 
($\bf A$,  $\bf O$, or $\bf S$) is
(canonically) {\em representable} 
if for every canonical
 matrix  realization of
$A_1$ there exists a canonical matrix realization of $A_2$ 
such that the matrix embedding $\ep$ is canonical. 
\end{ddd}

\begin{rrr} \label{invhom}
(1) It is easy to see that canonical matrix homomorphisms (\ref{e82})-(\ref{e7})
commute with the canonical matrix involutions (\ref{orth})-(\ref{atype})
(e.g. in type $\bf O$ the canonical involution is just matrix transpose).

(2) Note that compositions of canonical matrix homomorphisms are canonical. 
\end{rrr}

We are going to show that all embeddings
of involution simple algebras of the same type 
are representable, except for types $\bf O$ and $\bf S$ in characteristic 2.

\begin{ppp}\label{ema}
Let $\ep: A_1\to A_2$ be an embedding of 
finite dimensional involution simple algebras over $\bbF$ 
of type $\bf A$. Then $\ep$ is {\em representable}.
\end{ppp}

\begin{pf} Let $n_i$ be the degree of $A_i$.
Fix any bases of $V_1$ and $W_1$ such that the corresponding matrix realization
$M_{n_1}(\bbF)\oplus M_{n_1}(\bbF)$ of  $A_1$ is canonical
(i.e. the involution has the form (\ref{atype})).
Let $\pi_B$ (resp. $\pi_C$) denote the projection
$A_2\to B_2$ (resp. $A_2\to C_2$).
Fix any basis of $V_2$ which agree with the bases of 
$V_1$ and $W_1$ and the decomposition (\ref{sign}),
i.e. the projection $\pi_B\ep(A_1)$ has the following matrix form.
$$
\pi_B\ep(M,N)=
\diag(\underbrace{M,\dots,M}_{l},
\underbrace{N,\dots,N}_{r},\underbrace{0,\dots,0}_{z}),
\quad M,N\in M_{n_1}(\bbF).
$$
Fix a basis of $W_2$ such that the corresponding 
matrix realization of  $(A_2,*)$ is canonical 
(see Proposition~\ref{invol2}). Then
$$
\ep(M,N)=\ep((N^t,M^t)^*)=(\ep(N^t,M^t))^*=
((\pi_C\ep(N^t,M^t))^t,(\pi_B\ep(N^t,M^t))^t),
$$
so 
$$
\pi_C\ep(M,N)=(\pi_B\ep(N^t,M^t))^t=
\diag(\underbrace{N,\dots,N}_{l},
\underbrace{M,\dots,M}_{r},\underbrace{0,\dots,0}_{z}).
$$
Therefore
$$
\ep(M,N)=
(\diag(\underbrace{M,\dots,M}_{l},
\underbrace{N,\dots,N}_{r},\underbrace{0,\dots,0}_{z}),
\diag(\underbrace{N,\dots,N}_{l},
\underbrace{M,\dots,M}_{r},
\underbrace{0,\dots,0}_{z}))
$$
as required. 
\end{pf}

Our aim now is to prove a similar result
for orthogonal and symplectic
algebras in characteristic  $\ne2$. 
We need some auxiliary lemmas.

\begin{lll}[{\cite[2.23]{book}}]\label{tensor}
Let $D_1$ and $D_2$ be finite dimensional simple algebras
over $\bbF$ with involutions $\alpha_1$ and $\alpha_2$, respectively.
Then $\alpha=\alpha_1\otimes\alpha_2$ is an involution of
$D_1\otimes_\bbF D_2$. 
\begin{enumerate}
\item[$(i)$] 
If $\alpha_1$ and $\alpha_2$ are orthogonal,
then $\alpha$ is orthogonal.
\item[$(ii)$] 
If $\alpha_1$ is orthogonal and $\alpha_2$ is symplectic,
then $\alpha$ is symplectic.
\item[$(iii)$] 
If $\alpha_1$ and $\alpha_2$ are symplectic, 
then $\alpha$ is orthogonal in the case of $\chf\ne2$
and symplectic otherwise.
\end{enumerate}
\end{lll}

Recall that $e_1$ is the identity element of $A_1$.
We will use the notation $\bar A_1=e_1A_2e_1$ and $\bar V_1=e_1V_2$.
Let $\bar b_1$ be the restriction of the form $b_2$
to $\bar V_1$.

\begin{lll}\label{nond}
Assume that $A_2$ is not of type $\bf A$. Then 
\begin{enumerate}
\item[$(i)$] 
$\bar A_1$ is a $*$-invariant simple subalgebra of  $A_2$;
\item[$(ii)$] 
$\bar V_1$ is an irreducible $\bar A_1$-module and
$\bar V_1=\bar A_1 V_2$;
\item[$(iii)$] 
the form $\bar b_1$ on $\bar V_1$ is nondegenerate
and corresponds to the involution $*$ on $\bar A_1$;
moreover, 
$\bar b_1$ has the same type as $b_2$ except in the case when
$\chf=2$ and $A_2$ is of type $\bf O$.
\end{enumerate}
\end{lll}

\begin{pf} Note that $e_1$ is an idempotent of $A_2$
and $e_1^*=e_1$, so $(i)$ and $(ii)$ are clear. Now assume
that $\bar b_1$ is degenerate, i.e. there exists
$v\in V_2$ such that $e_1v\ne0$ and 
$b_2(e_1v,e_1w)=0$ for all $w\in V_2$. Then
$$
b_2(e_1v,w)=b_2(e_1e_1v,w)=b_2(e_1v,e_1w)=0 \quad{\rm for\ all\ } 
w\in V_2,
$$
which contradicts to nondegeneracy of $b_2$.
It remains to note that if $b_2$ is
alternating (resp. symmetric), then $\bar b_1$
is
alternating (resp. symmetric). 
\end{pf}

\begin{lll}\label{deco}
Let $\ep:A_1\to A_2$ be an embedding of involution simple
algebras of types different from $\bf A$
and let $\alpha_i$ denotes the involution of $A_i$. 
Fix any canonical matrix realization 
$(M_{n_1}(\bbF),\alpha_1)$ of $A_1$. Then there exists
a matrix realization 
$(M_{n_2}(\bbF),\alpha_2)$ of $A_2$ such that 
the following hold.
\begin{enumerate}
\item[$(i)$] 
The embedding $\ep$ is
the composition of the following embeddings
of algebras with involution.
$$
(M_{n_1}(\bbF),\alpha_1)\stackrel{\eta}{\longrightarrow}
(M_{n_1}(\bbF)\otimes_\bbF M_k(\bbF),\alpha_1\otimes\beta_1)
\stackrel{\iota}{\longrightarrow}
(M_{kn_1}(\bbF),\beta_2)\stackrel{\zeta}{\longrightarrow}
(M_{n_2}(\bbF),\alpha_2)
$$
where $\eta(X)=X\otimes e$
with $e$ the identity element of $M_k(\bbF)$,
$\iota$ is the natural isomorphism, and $\zeta$
is a natural embedding (i.e. of signature $(1,0,z)$).
\item[$(ii)$] 
If $\chf\ne2$ and $A_1$ and $A_2$ 
are both of type $\bf O$, then $\alpha_1=\beta_1=\beta_2=\alpha_2=t$ (matrix
transpose)
\item[$(iii)$] 
If $\chf\ne2$ and $A_1$ and $A_2$ 
are both of type $\bf S$, then $\alpha_1=\beta_2=\alpha_2=\tau$ 
(symplectic transpose) and 
$\beta_1=t$.
\end{enumerate}
\end{lll}

\begin{pf}
Let $\bar A_1$ be as in Lemma~\ref{nond} and let 
$\beta_2$ be the restriction of $\alpha_2$ to $\bar A_1$. Let
$B$ be the centralizer
of $A_1$ in $\bar A_1$. Note that  $A_1$ and $\bar A_1$ 
are simple and have the same identity element. Therefore
$B$ is simple and 
$\bar A_1=A_1B\cong A_1\otimes_\bbF B$
 (see e.g. \cite[1.5]{book}).
Clearly, $B$ is $\beta_2$-invariant. 
Denote by $\beta_1$ the restriction of $\beta_2$ to $B$.
Then 
$(\bar A_1,\beta_2)\cong (A_1\otimes_\bbF B,\alpha_1\otimes\alpha_2)$.
We get the following chain of embeddings of
algebras with involution:
$$
A_1\longrightarrow A_1\otimes_\bbF B\simeq \bar A_1\longrightarrow A_2.
$$
Identifying $A_1$ with $M_{n_1}(\bbF)$,
$B$ with $M_k(\bbF)$ for some $k$, $\bar A_1$
with $M_{kn_1}(\bbF)$, and $A_2$ with $M_{n_2}(\bbF)$,
we prove (i).

Assume now that $\chf\ne2$  and $A_1$ and $A_2$ 
are of the same type $\bf S$ (resp. $\bf O$), i.e.
$\alpha_1$ and $\alpha_2$ are of type $\bf S$ (resp. $\bf O$).
Then by Lemma~\ref{nond}, $\beta_2$
is of type $\bf S$ (resp. $\bf O$). Therefore by Lemma~\ref{tensor},
$\beta_1$ is of type $\bf O$.
Fixing an appropriate isomorphism $B\cong M_k(\bbF)$,
by Lemma~\ref{invol1}, we can assume that $\beta_1$
is a matrix transpose. Using the same lemma we get
that $\beta_2$ can be represented as $\tau$ (resp. $t$).
Now by Lemma~\ref{nond}, the restriction of the form $b_1$ 
to $\bar V_1$ is nondegenerate.
Thus $V_2=\bar V_1\oplus \bar V_1^\perp$. By choosing 
a suitable basis in $\bar V_1^\perp$, we can easily
represent $\alpha_2$ as $\tau$ (resp. $t$).
\end{pf}

As a corollary we get the following 
analogue of Proposition~\ref{ema}
for symplectic and orthogonal algebras.

\begin{ppp}\label{emb}
Let $\ep: A_1\to A_2$ be an embedding of 
finite dimensional involution simple algebras over $\bbF$ 
of the same type $\bf S$ or $\bf O$. Assume that $\chf\ne2$. 
Then $\ep$ is {\em representable}.
That is, for every canonical
 matrix  realization of
$A_1$ there exists a canonical matrix realization of $A_2$ 
such that the embedding $\ep$ is of the form (\ref{e82}).
\end{ppp}

Proposition \ref{ch2} below shows that the case of characteristic $2$ is
exceptional indeed.

We will also need the following result, which describes
embeddings
of involution simple algebras of different types.
Recall that $(l,r,z)$
is the signature of the embedding $\ep:A_1\to A_2$.

\begin{ppp}\label{p3}
Let $\ep:A_1\to A_2$ be an embedding of 
finite dimensional involution simple algebras over $\bbF$. Assume that $\chf\ne2$.
\begin{enumerate}
\item[(i)] If $A_1$ is of type $\bf A$ and $A_2$ is not of type $\bf A$, 
then $l=r$.
\item[(ii)]If $A_1$ is of type $\bf S$ (resp., $\bf O$) and $A_2$ is 
of type $\bf O$ (resp., $\bf S$), then $l$ is even.
\item[(iii)] If $A_1$ and $A_2$ are both not of type $\bf A$ and $l$ 
is even, then there exist an algebra $D$ of type $\bf A$,
an embedding $\eta:A_1\to D$ with the signature $(l/2,0,0)$ and an embedding
$\zeta:D\to A_2$ with the signature $(1,1,z)$ such that 
$\ep=\zeta\eta$.
\item[(iv)] If $A_1$ and $A_2$ are of type $\bf A$ and $l=r$, 
then there 
exist an algebra $D$ of type $\bf O$ (resp., $\bf S$),
embeddings $\eta:A_1\to D$ of signature $(l,l,0)$ and
$\zeta:D\to A_2$ of signature $(1,0,z)$ such that 
$\ep=\zeta\eta$.
\end{enumerate}
\end{ppp}

\begin{pf}
$(i)$
Recall that $A_i=B_i\oplus C_i$ where $B_i$ and $C_i$ are the simple components
of $A_i$ and $B_i^*=C_i$. 
And $f_i$ and $g_i=f_i^*$ are the identities of $B_i$ and $C_i$, respectively.
Obviously, $l=(\dim f_1A_2f_1)/n_1$
and $r=(\dim g_1A_2g_1)/n_1$ where $n_1=\dim V_1=\dim W_1$.
Since $(f_1A_2f_1)*=g_1A_2g_1$, we get that $l=r$.
Note that this is valid for the case of $\chf=2$ as well.

$(ii)$ Represent the embedding $A_1\to A_2$
as in Lemma~\ref{deco}$(i)$. Note that $k=l$.
By Lemma~\ref{nond}$(iii)$, $\beta_2$ has the same type as 
$\alpha_2$. Thus the types of $\alpha_1$ and $\beta_2$ are different.
By Lemma~\ref{tensor}, $\beta_1$ must be symplectic.
Therefore $k=l$ is even.

$(iii)$ Represent the embedding $A_1\to A_2$
as in Lemma~\ref{deco}$(i)$. Denote by $B$ the algebra $M_k(\bbF)$.
By assumption, $k=l$ is even.
Let $e$ be the identity element of $B$. By Lemma~\ref{fg},
$B$ has two idempotents 
$f$ and $g$ such that $e=f+g$, $fg=gf=0$, and $f^*=g$.
Then $B_f=fBf$ and $B_g=gBg$ are simple subalgebras of $B$,
$B_f\cap B_g=0$, $B_fB_g=B_gB_f=0$, and $B_f^*=B_g$. Thus
$B'=B_f\oplus B_g$ is an involution simple subalgebra of $B$
of type $\bf A$. Therefore $D=A_1\otimes_\bbF B'$ is
an involution simple subalgebra of $A_1\otimes_\bbF B$
of type $\bf A$. Since $e=f+g$, $D$ contains $A_1$.
Clearly the signature of the embedding
$A_1\to D$ is $(l/2,0,0)$ and the signature of the embedding
$D\to A_2$ is $(1,1,z)$.

$(iv)$ Let $A=M_n(\bbF)\oplus M_n(\bbF)$
be an involution simple algebra with standard involution
$(X,Y)^*=(Y^t,X^t)$. Let $k\le n$ and let the algebra $D=M_k(\bbF)$
have an involution $\alpha$. Define a "corner"
embedding $\varphi:D\to A$ via $\varphi(Z)=(\bar Z, \overline{(Z^\alpha)^t})$
where $\bar Z=\diag(Z,0\dots,0)$. Since $t\circ \alpha$ is an 
automorphism of $D$, $\varphi$ is an algebra homomorphism.
Moreover, one can easily check that $\varphi$ respects 
involution:
$$
\varphi(Z^\alpha)=(\overline{Z^\alpha}, \overline{Z^t})=
(\bar Z, \overline{(Z^\alpha)^t})^*=\varphi(Z)^*
$$
By Proposition~\ref{ema}, the embedding $\ep$
can be represented as in (\ref{e7}) with $l=r$ and involution $*$
acting as $(X,Y)^*=(Y^t,X^t)$ on both algebras. 
Now let $\alpha$ be either symplectic involution $\theta_-$
or orthogonal involution $\theta_+$ (see Remark~\ref{other})
of the algebra $D=M_{2l}(\bbF)$. Let $\zeta=\varphi:D\to A_2$
be the corner embedding of algebras with involution described above.
Note that it is an embedding of signature $(1,0,z)$.
Observe that 
$$
\left(\left(
\begin{array}{cc}
a & b \\
c & d
\end{array}
\right)^{\theta_\pm}\right)^t=
\left(
\begin{array}{cc}
d^t & \pm b^t \\
\pm c^t & a^t
\end{array}
\right)^t=
\left(
\begin{array}{cc}
d & \pm c \\
\pm b & a
\end{array}
\right),
\quad {\rm for \ } 
\left(
\begin{array}{cc}
a & b \\
c & d
\end{array}
\right)\in M_{2l}(\bbF),
$$
where $a,b,c,d$ are square matrices of size $l$. 
Therefore, it is easy to see from formula (\ref{e7})
that $\ep(A_1)\subset \zeta(D)$. Define by $\eta$
the following composition of embeddings:
$$
A_1\longrightarrow \ep(A_1) \longrightarrow 
\varphi(D) \stackrel{\zeta^{-1}}{\longrightarrow}D.
$$
Then $\eta$ is of signature $(l,r,0)$ and $\ep=\zeta\eta$, as required.
\end{pf}

It remains to consider the case of characteristic 2, which is a bit
more complicated. 

\begin{lll}\label{iden}
 Let $\chf=2$ and
let $\ep:A_1\to A_2$ be an
embedding of involution simple algebras preserving the identity element
(i.e. $\ep(e_1)=e_2$). Assume that
$A_2$ is of type $\bf O$. Then $A_1$ is of type $\bf O$.
\end{lll}

\begin{pf} Assume that $A_1$ is of type $\bf A$  or
$\bf S$. Then by Lemma~\ref{fg}, $A_1$ has idempotents $f$ and $g$ such that $e_1=f+g$, 
$fg=gf=0$, and $f^*=g$.  Let $b$ be a symmetric nondegenerate
form on $V_2$ corresponding to the involution.
Then for all $v\in V_2$
we have 
$$
b(v,v)=b((f+g)v,v)=b(fv,v)+b(v,fv)=0,
$$
as $b$ is symmetric. Therefore $b$ is alternating, so $A_2$
is symplectic, which contradicts the assumption.
\end{pf}

\begin{ppp}\label{ch2} Let $\chf=2$ and
let $\ep:A_1\to A_2$ be an
embedding of involution simple algebras of the same type $X=\bf O$ or 
$\bf S$.
Then the following conditions are equivalent.
\begin{enumerate}
\item[(i)] The embedding $\ep$ is representable.
\item[(ii)]  Each $*$-invariant involution
simple subalgebra $D$ of $A_2$ containing  $A_1$
is of type $X$.
\end{enumerate}
Moreover, if the embedding $\ep$ is not representable, then
there exists a $*$-invariant involution
simple subalgebra $D$ of $A_2$ which is 
of type $\bf A$
and contains $A_1$.
\end{ppp}

\begin{pf} 
By Lemma~\ref{deco}$(i)$, the embedding $\ep$ can be
represented as the composition of embeddings
$A_1\to C \to A_2$ where $C=e_1A_2e_1\cong A_1\otimes M_k(\bbF)$ 
is involution simple 
of type $\bf O$ or $\bf S$ and
has the same identity  element $e_1$ 
as $A_1$, and the embedding $C\to A_2$ is natural
(of signature $(1,0,z)$).

$(i)\Rightarrow(ii)$ ($X=\bf O$):
Assume that  $\ep$ is representable and
there exists a 
$*$-invariant involution
simple subalgebra $D$ of $A_2$ containing  $A_1$ which 
is not of type $\bf O$.
The matrix presentation (\ref{e82}) shows that
$C$ is of type $\bf O$.
Let $e_D$ be the identity  element of $D$.
Since $e_D$ is an idempotent, the algebra
$F=e_DA_2e_D$ is a $*$-invariant simple
subalgebra of $A_2$ containing $D$. By Lemma~\ref{iden},
it cannot be orthogonal. Therefore $F$ is of type
$\bf S$. Note that $F$ contains $C$
and $e_1Fe_1=C$. Therefore by Lemma~\ref{nond}$(iii)$,
$C$ must be of the same type $\bf S$, which is a contradiction.

$(i)\Rightarrow(ii)$ ($X=\bf S$):
Assume that  $\ep$ is representable and
there exists a 
$*$-invariant involution
simple subalgebra $D$ of $A_2$ containing  $A_1$ which 
is not of type $\bf S$. 
First assume that $D$ is of type $\bf A$.
Then $e_1De_1$ is an involution simple subalgebra 
of $C=e_1A_2e_1$
of type $\bf A$ containing $A_1$.
Recall that $C\cong A_1\otimes M_k(\bbF)$.
Therefore $e_1De_1\cong A_1\otimes E$ where
$E$ is an involution simple subalgebra of $M_k(\bbF)$
of type $\bf A$ with the same identity  element.
Since $\ep$ is representable, the involution on $M_k(\bbF)$
is orthogonal, which contradicts to Lemma~\ref{iden}.

Suppose now that $D$ is of type $\bf O$.
As in the case $ X=\bf O$, 
the algebra $F=e_DA_2e_D$ is a $*$-invariant simple
subalgebra of $A_2$ containing $C$.
By Lemma~\ref{nond}$(iii)$, $F$ is symplectic.
Therefore $F\cong D\otimes_\bbF M_q(\bbF)$
with a symplectic involution on $M_q(\bbF)$ (Lemma \ref{tensor}(i)). 
By Lemma~\ref{fg}, $M_q(\bbF)$ has two idempotents
$f$ and $g$ such that $f+g$ is the identity element of $M_q(\bbF)$,
$fg=gf=0$, and $f^*=g$. Therefore 
$D'=D\otimes f \oplus D\otimes g$ is an involution simple
subalgebra of $A_2$ of type $\bf A$ containing
$A_1$. However the case of type $\bf A$ subalgebra containing $A_1$ has been 
already considered in the previous paragraph.

$(ii)\Rightarrow(i)$ and "Moreover" part: Assume that 
the embedding $\ep$ is not representable. 
We are going to show that $A_2$ contains an involution simple
subalgebra of type $\bf A$ containing
$A_1$. 
Recall that $C\cong A_1\otimes M_k(\bbF)$
and the restriction of the involution 
$*$ on $C$ has the form
$\alpha_1\otimes \alpha_2$ where $\alpha_1$ is
the involution of $A_1$ and $\alpha_2$
is an involution of $M_k(\bbF)$. 
Clearly if 
$\alpha_2$ is orthogonal, then $\ep$ is representable
(see the proof of Lemma~\ref{deco}). 
Therefore $\alpha_2$  is symplectic.
Then, as above, $M_k(\bbF)$ has two idempotents
$f$ and $g$ such that $f+g$ is the identity element of $M_k(\bbF)$,
$fg=gf=0$, and $f^*=g$. Therefore 
$D=A_1\otimes f \oplus A_1\otimes g$ is an involution simple
subalgebra of $A_2$ of type $\bf A$ containing
$A_1$. The proposition follows.
\end{pf}

The following results show how embedding signatures behave under compositions.

\begin{ppp}\label{p70}
Let $\ep_1: A_1\to  A_2$ and $\ep_2: A_2\to  A_3$ be  
embeddings of involution simple algebras of the same type with the signatures 
$(l_1,r_1,z_1)$ and $(l_2,r_2,z_2)$, respectively. 
Denote by $(l,r,z)$ the 
signature of $\ep=\ep_2\ep_1$. Then
\begin{eqnarray}
l&=&l_1 l_2+r_1 r_2, \label{e88}\\
r&=&r_1 l_2+l_1 r_2, \label{e89}\\
z&=&z_1(l_2+r_2)+z_2.\nonumber
\end{eqnarray}
\end{ppp}

\begin{pf}
For types $\bf S$ and $\bf O$ one has $r=r_1=r_2=0$, so the
statement immediately follows from (\ref{sign}). For type
$\bf A$, the embeddings are representable so one can use
(\ref{e7}).
\end{pf}

Note that $l+r=(l_1+r_1)(l_2+r_2)$ and $l-r=(l_1-r_1)(l_2-r_2)$.
Thus, the following is true.

\begin{ccc}\label{c80}
Let $ A_1\to\dots\to  A_k$ be a sequence 
of embeddings of involution simple algebras of the same type.
Let $(l_i,r_i,z_i)$ be the signature of $ A_i\to A_{i+1}$, $(l,r,z)$ the 
signature of $ A_1\to  A_k$,
$s_i=l_i+r_i$, $c_i=l_i-r_i$, $s=l+r$, $c=l-r$. Then
$s=s_1\dots s_{k-1}$ and $c=c_1\dots c_{k-1}$.
\end{ccc}

Recall that $n_i$ is the degree of $A_i$ (so 
$A_i\cong M_{n_i}(\bbF)$ or $M_{n_i}(\bbF)\oplus M_{n_i}(\bbF)$).

\begin{lll}\label{l25}
Let $\ep_1: A_1\to  A_2$ and $\ep: A_1\to  A_3$ be 
representable embeddings of involution simple algebras of the same type 
($\bf A$, $\bf S$ or $\bf O$) with the 
signatures $(l_1,r_1,z_1)$
and $(l,r,z)$, respectively. Assume that a triple of non-negative integers 
$(l_2,r_2,z_2)$ satisfies the following conditions
\begin{eqnarray}
l+r&=&(l_1+r_1)(l_2+r_2), \label{e40}\\
l-r&=&(l_1-r_1)(l_2-r_2), \label{e41}\\
n_3&=&n_2(l_2+r_2)+z_2  \label{n11}
\end{eqnarray}
where $n_i$ is the degree of $A_i$. 
Then there exists a representable embedding 
$\ep_2: A_2\to A_3$ with the signature $(l_2,r_2,z_2)$ such that 
$\ep=\ep_2\ep_1$.
\end{lll}

\begin{pf}
Fix any canonical matrix realizations of $A_1$, $A_2$, $A_3$ such that
the matrix embeddings $\ep_1$ and $\ep$ become canonical (see Definition \ref{canembed}). Consider the canonical matrix embedding $\ep_2:A_2\to A_3$
with signature $(l_2,r_2,z_2)$.  The embedding $\ep_2$ is well-defined because of 
(\ref{n11}) and respects the involution (see Remark \ref{invhom}(1)).
By Remark \ref{invhom}(2), the matrix homomorphism $\ep_2\ep_1$ is canonical.
By rewriting the conditions (\ref{e40}) and (\ref{e41})  in the form 
(\ref{e88}) and (\ref{e89}) we see  that both canonical
homomorphisms $\ep$ and $\ep_2\ep_1$ have
the same signature, so $\ep=\ep_2\ep_1$.
\end{pf}

Our classification will be given in terms of so-called supernatural (or Steinitz)
numbers. They are defined as follows.
Let $(p_1,p_2,\dots)$ be the increasing sequence of all prime numbers. The 
set of all mappings from $\{p_1,p_2,\dots\}$ into the set 
$\{0,1,2,\dots\}\cup\{\infty\}$ is called the set of {\em supernatural} 
(or {\em Steinitz}) numbers. 
If a supernatural number takes a value $\alpha_1$ at $p_1$, $\alpha_2$ at 
$p_2$,\dots, this element will be denoted by 
$p_1^{\alpha_1}p_2^{\alpha_2}\dots$. 
The supernatural numbers can be regarded as formal products of powers
of primes where infinity
is permitted as a power. 
The set of natural numbers  
$\bbN$ can be identified in an evident way with a subset of supernatural 
numbers. If $\Pi=p_1^{\alpha_1}p_2^{\alpha_2}\dots$ and 
$\Pi'=p_1^{\alpha_1'}p_2^{\alpha_2'}\dots$ are two supernatural numbers, we 
set $\Pi\Pi'=p_1^{\alpha_1+\alpha_1'}p_2^{\alpha_2+\alpha_2'}\dots$.
We say that $\Pi$ divides $\Pi'$ if and only if 
$\alpha_1\le\alpha_1',\alpha_2\le\alpha_2'\dots$. 
Let $q\in\bbQ$. We write $\Pi=q\Pi'$ (or $q\in\frac{\Pi}{\Pi'}$) if there 
exists $n\in\bbN$ such that $nq\in\bbN$ and $n\Pi=nq\Pi'$. 
If there exists non-zero $q\in\bbQ$ such that $\Pi=q\Pi'$, then we say that 
$\Pi$ and $\Pi'$ are {\em $\bbQ$-equivalent} and denote this relation by 
$\Pi\stackrel{\bbQ}{\sim}\Pi'$.
Let $\Se=(s_1,s_2,\dots)$ be a sequence of natural numbers. Denote by 
$\Pi(\Se)$ the supernatural number $s_1s_2s_3\dots$.
We will use the following simple observation.

\begin{ppp}[{\cite[Proposition 3.2]{bzh}}]
\label{p33}
Let $\Se=(s_i)_{i\in I}$ and $\Se'=(s_j')_{j\in J}$ be sequences of natural numbers. 
Then $q\in\frac{\Pi(\Se)}{\Pi(\Se')}$ if and only if for each $i\in I$ and 
$k\in J$ there exist $j=j(i)\in J$ and $l=l(k)\in I$ such that $s_1\dots 
s_i$ divides $q s_1'\dots s_j'$ (over $\Ze$) and $q s_1'\dots s_k'$ divides 
$s_1\dots s_l$ (over $\Ze$).
\end{ppp} 

%\begin{pf}
%Let $q\in\frac{\Pi(\Se)}{\Pi(\Se')}$. This is equivalent to say that 
%there exists $n\in\bbN$ such that $nq\in\bbN$ and $n\Pi(\Se)=nq\Pi(\Se')$. 
%It is not difficult to see that the latter property is equivalent to the 
%following. For each $i\in I$ and $k\in J$ there exist $j=j(i)\in J$ and 
%$l=l(k)\in I$ such that $ns_1\dots s_i$ divides $nq s_1'\dots s_j'$ (over 
%$\Ze$) and $nq s_1'\dots s_k'$ divides $ns_1\dots s_l$ (over $\Ze$). So the 
%lemma follows.
%\end{pf}

\section{Bratteli diagrams and dimension groups}
\label{dirlimits}

Let 
\begin{equation} \label{e000}
A_1\to A_2\to A_3\to\dots\to A_i\to A_{i+1}\to\dots
\end{equation}
be a sequence of embeddings of finite dimensional involution simple algebras over $\bbF$. 
Assume that all $A_i$ are of the same type and $\chf\ne2$. 
Then, as we proved in Propositions~\ref{ema} and \ref{emb},
all embeddings $A_i\to A_{i+1}$ are representable, so one can assume that all $A_i$ are
  matrix (or double matrix) algebras and the embeddings and involutions are canonical.  
This justifies the following definition.

\begin{ddd}\label{limtype} 
Let $A$ be a locally involution simple associative algebra of countable dimension. We say that
$A$ is {\em canonically representable} if 
it is isomorphic to the direct limit of the sequence (\ref{e000}) where all 
$A_i$ are matrix (resp. double matrix) algebras with canonical involutions
of the same type $X$ ($={\bf A},{\bf S}$ or $\bf O$) and all embeddings are canonical.
In that case we say that the sequence (\ref{e000}) 
is a {\em canonical representation} for $A$ and 
$A$ is of {\em type} $X$.
\end{ddd}

Note that the type $X$ of the algebra $A$ may not be unique. 

Proposition~\ref{ch2} shows that some of the embeddings $A_i\to A_{i+1}$ may not be 
representable  in characteristic 2. 
Fortunately, there is a way to modify the sequence
(\ref{e000}), without changing the limit algebra, in order to 
get representable embeddings  even in characteristic 2.

\begin{ttt}\label{rt}
Let $A$ be a locally involution simple associative algebra over $\bbF$ of countable dimension. 
Then $A$ is canonically representable.
\end{ttt}

\begin{pf}
Let 
$A_1\to A_2 \to A_3\to \dots$ be a sequence 
of embeddings 
of involution simple finite dimensional 
associative algebras such that $A=\dlim A_i$. 
Choose an infinite subsequence of algebras of the same type.
If $\chf\ne2$, or $\chf=2$ and all algebras are of type $\bf A$,
then all embeddings are representable by Propositions \ref{ema} and \ref{emb}. 
Assume $\chf=2$. If there is an infinite number of non-representable
embeddings, then by Proposition~\ref{ch2}, we can
replace the subsequence by a sequence of embeddings 
of algebras of type $\bf A$. 
Otherwise, we get the result by removing a finite number of algebras
in the beginning of the sequence.
\end{pf}

Theorem~\ref{rt} reduces classification of locally 
involution simple algebras to the following two problems:
\begin{itemize}
\item[(a)] classification of the direct limits of canonical sequences 
of the same type; 
\item[(b)] classification of intertype isomorphisms.
\end{itemize}

We are going to simplify Problem (a) even further and reduce it to the algebras without involution. 
We need the following trivial observation. 

\begin{ppp}\label{p50}  Let 
$A_1\to A_2 \to A_3\to \dots$ and 
$A_1'\to A_2' \to A_3'\to \dots$
be two sequences of embeddings 
of algebras (or algebras with involution). Then
$\dlim A_i\cong\dlim A_j'$ if and only if 
there exist sequences of indices $i_1<i_2<\dots$ and
$j_1<j_2<\dots$ and homomorphisms $\varphi_k: A_{i_k}\to A_{j_k}'$
and
$\varphi_k': A_{j_k}'\to A_{i_{k+1}}$ ($k=1,2,\dots$) 
such that the following diagram commutes.
\begin{equation}\label{iso}
\begin{array}{cccccccccc}
 A_{i_1} & \longrightarrow &  A_{i_2} & \longrightarrow & \dots &
 A_{i_k} & \longrightarrow &  A_{i_{k+1}} & \longrightarrow & \dots \\
{\downarrow}\lefteqn{\scriptstyle \varphi_1} &
{\nearrow} \lefteqn{\scriptstyle \varphi_1'}&
{\downarrow} \lefteqn{\scriptstyle \varphi_2}&
{\nearrow} \lefteqn{\scriptstyle \varphi_2'}&  &
{\downarrow}\lefteqn{\scriptstyle \varphi_k} &
{\nearrow} \lefteqn{\scriptstyle \varphi_k'}&
{\downarrow} \lefteqn{\scriptstyle \varphi_{k+1}}&
{\nearrow} \lefteqn{\scriptstyle \varphi_{k+1}'}& \\
 A_{j_1}' & \longrightarrow &  A_{j_2}' & \longrightarrow & \dots &
 A_{j_k}' & \longrightarrow &  A_{j_{k+1}}' & \longrightarrow & \dots 
\end{array}
\end{equation}
\end{ppp}

\begin{pf}
Set $ A=\dlim A_i$ and $ A'=\dlim A_j'$. Assume that there exists an 
isomorphism $\varphi: A\to A'$. Fix any index $i_1$. Then there exists 
$j_1$ such that $\varphi( A_{i_1})\subseteq A_{j_1}'$. Similarly, there 
exists $i_2$ such that $\varphi^{-1}( A_{j_1}')\subseteq A_{i_2}$, and 
so on. 
Denote by $\varphi_k$ the restriction of $\varphi$ to $A_{i_k}$, 
and by $\varphi_k'$ the restriction of $\varphi^{-1}$ to $A_{j_k}'$, 
$k=1,2,\dots$. Then the diagram above commutes.
The converse statement is obvious.
\end{pf}

\begin{ttt} \label{aareduced}
Two locally involution simple associative algebras of the same type over $\bbF$ of countable dimension
are isomorphic if and only if they are isomorphic as associative algebras.
\end{ttt}

\begin{pf} 
Let $ A$ and $ A'$ be two locally involution simple associative algebras
and let  $ A=\dlim A_i$ and $ A'=\dlim A_j'$ be their canonical 
representations. 
Assume that  $ A$ and $ A'$ are isomorphic as associative algebras. Using Proposition \ref{p50}, 
we get a commutative diagram (\ref{iso}), where $\varphi_k:A_{i_k}\to A_{j_k}'$
and $\varphi_k':A_{j_k}'\to A_{i_{k+1}}$ are algebra homomorphisms, 
not necessarily respecting the involution. 
Let $\ep_k: A_{i_k}\to A_{i_{k+1}}$ and $\ep_k': A_{j_k}'\to A_{j_{k+1}}'$
be the horizontal maps. Note that they are canonical and respect the involution. 
Denote
by $\psi_k$ (resp. $\psi_k'$) the canonical map 
$A_{i_k}\to A_{j_k}'$ (resp. $A_{j_k}'\to A_{i_{k+1}}$)
of the same signature as $\varphi_k$ (resp. $\varphi_k'$). 
Then by Remark \ref{invhom}(1) these maps respect the involution.
It remains to show that they make the diagram (\ref{iso}) commutative. 
Note that the signature of $\psi_k'\psi_k$ equals to the signature of 
$\varphi_k'\varphi_k=\ep_k$. 
Since both $\psi_k'\psi_k$ and $\ep_k$ are canonical, we get that $\psi_k'\psi_k=\ep_k$.
Similarly, one proves that $\psi_{k+1}\psi_k'=\ep_k'$. Therefore 
the diagram (\ref{iso}) commutes with respect to the maps
$\psi_k$ and $\psi_k'$. Proposition \ref{p50} implies 
that $ A$ and $ A'$ are isomorphic as algebras with involution.

The converse statement is trivial.
\end{pf}

Theorem \ref{aareduced} reduces Problem (a) to classifying the direct limits 
of finite dimensional semisimple algebras. This is usually done
in terms of Bratteli diagrams, the $K_0$-functor and dimension groups.
To make the statements of the results a little bit clearer  it is best
to work in the category of unital algebras (i.e. algebras with identity elements and with  identity preserving  homomorphisms). In our case 
this can be easily achieved by adjoining an external identity element. 

\begin{ddd}\label{extid}
Let $A$ be an associative algebra. Define the algebra $\hat{A}$ as follows.
 If $A$ has an identity element, put $\hat{A}=A$.
Otherwise, put $\hat{A}=A+ \bbF {\bf 1}_{\hat{A}}$ where 
${\bf 1}_{\hat{A}}$ is the identity element of $\hat{A}$.
\end{ddd}

Note that if $A$ has an involution then this involution trivially extends to $\hat{A}$.

\begin{lll}\label{unital}
Let $A$ be a locally semisimple associative algebra.
Then $\hat{A}$ is locally semisimple in the category of unital algebras.
\end{lll}

\begin{pf} Let $A=\dlim A_i$ with $A_i$ finite dimensional semisimple.
If $A$  contains an identity element ${\bf 1}_A$, then $A$ is the direct limit of  those
$A_i$ which contain ${\bf 1}_A$, as required.
If $A$  has no identity element, then 
$\hat{A}=A+ \bbF {\bf 1}_{\hat{A}}=\dlim B_i$
where $B_i=A_i+\bbF {\bf 1}_{\hat{A}}$ are obviously finite dimensional and semisimple.
\end{pf}

\begin{ppp} \label{adiso}
Let $A$ and $A'$ be involution simple associative algebras.
Then $A\cong A'$ as associative algebras if and only if
$\hat{A}\cong \hat{A}'$ as associative algebras. 
\end{ppp}

\begin{pf} By construction, $A\cong A'$ implies 
$\hat{A}\cong \hat{A}'$. 
Assume now that $\hat{A}\cong \hat{A}'$. We need to show that $A\cong A'$.
Denote by $\Soc(A)$ the sum of all minimal ideals of $A$. 
Then by Proposition \ref{issum}, $\Soc(A)=A$.
Obviously, $\Soc(A)\subseteq \Soc(\hat{A})$.
We claim that  $\Soc(A)= \Soc(\hat{A})$. Indeed, this is obvious if $A=\hat{A}$.
Assume $A\ne\hat{A}$, i.e. $A$ has no identity element. 
Let $M$ be a minimal ideal of $\hat{A}$ such that
$M\not\subseteq \Soc(A)$. Then $M\cap\Soc(A)=0$.
But $\Soc(A)=A$ is an ideal of codimension $1$ in  $\hat{A}$.
Therefore $M$ is one-dimensional and $\hat{A}=A\oplus M $.
Write ${\bf 1}_{\hat{A}}=a+m$ where $a\in A$ and $m\in M$. 
Then obviously $a$ is an identity element of $A$, which is a contradiction.
Therefore, $A=\Soc(\hat{A})\cong\Soc(\hat{A'})=A'$, as required.
\end{pf}

Locally semisimple algebras are best described in terms of their Bratteli diagrams.
These are defined as follows.
Let $B$ be the direct limit of the infinite sequence 
\begin{equation}\label{bseq}
B_1\to B_2\to B_3 \to \dots
\end{equation}
where the $B_i$ are finite dimensional semisimple algebras over $\bbF$. 
Let $S_i^1,S_i^2,\dots S_i^{k_i}$ be the simple components of $B_i$, i.e. 
$B_i=S_i^1\oplus S_i^2\oplus \dots\oplus S_i^{k_i}$. Let $V_i^j$ be the natural
$S_i^j$-module. Then $V_i^j$ can be considered as an $B_i$-module. 
Denote by $m_i^{jq}$ the multiplicity of $V_i^j$ in the restriction
of $V_{i+1}^q$ to $S_i^j$ (i.e. $m_i^{jq}$ is the number of copies of
$S_i^j$ that are mapped to $S_{i+1}^q$). The Bratteli diagram
of the sequence (\ref{bseq}) 
consists of the vertices $V=\{V_i^j\mid i=1,2,3,\dots;\ 1\le j\le k_i\}$
and edges. Two vertices $V_i^j$ and $V_{i+1}^q$ are connected by an edge
if and only if $m_i^{jq}>0$. In that case the edge is labelled by the
number $m_i^{jq}$. Let $n_i^j=\dim V_i^j$ be the degree of $S_i^j$. Then obviously
\begin{equation}\label{mijq}
\sum_{j=1}^{k_i}m_i^{jq}n_i^j\le n_{i+1}^q
\end{equation}
Moreover, if all homomorphisms in (\ref{bseq}) 
are unital then we have equality in (\ref{mijq}) for all $i$ and $q$, so
the whole sequence (\ref{bseq})
can be reconstructed from its Bratteli diagram 
provided the degrees of the simple
components of the first term $B_1$ are known 
(in the case of non-unital embeddings extra data is needed). 

Now let $A$ be a locally involution simple associative algebra of type 
$X$ ($={\bf A},{\bf S}$ or $\bf O$)  
over $\bbF$ of countable dimension. 
By Theorem \ref{rt}, $A$ is the direct limit of the sequence 
(\ref{e000}) where all 
$A_i$ are matrix (resp. double matrix) algebras with canonical involutions
of the same type $X$ and all embeddings are canonical.

We will denote by $(l_i,r_i,z_i)$ 
 the signature of the embedding $ A_i\to A_{i+1}$ and by 
$n_i$ the degree of $A_i$ (i.e.  $A_i=M_{n_i}(\bbF)$ and $r_i=0$ for $X={\bf S},{\bf O}$
and $A_i=M_{n_i}(\bbF)\oplus M_{n_i}(\bbF) $ for $X={\bf A}$). 
By Remark~\ref{lr}, for type $\bf A$ algebras
we can and will assume  
that $l_i\ge r_i$ for all 
$i$. 
It is convenient to add to the sequence an algebra of degree 1
(the 1-dimensional algebra $\bbF$ is considered to be of both types
$\bf O$ and ${\bf S}$), 
so we will
assume that $n_1=1$, $l_1=n_2$ and $r_1=z_1=0$. 
Denote by $\Te$ the triple sequence $(l_i,r_i,z_i)_{i\in \bbN}$. 
Since $n_{i+1}=(l_i+r_i)n_i+z_i$ for all $i$, the canonical sequence (\ref{e000}) is uniquely determined by the triple sequence $\Te$ and type $X$. 
We will denote by $A(\Te,X)$ the corresponding locally involution simple
associative algebra over  $\bbF$, by 
$A(\Te)$ the corresponding locally semisimple algebra 
(i.e. the direct limit of the associative algebras (\ref{e000}) 
disregarding the involution) and by $\hat{A}(\Te)$ the corresponding algebra 
with an identity element (see Definition \ref{extid}). 
 Recall that by Theorem \ref{aareduced} and Proposition \ref{adiso},
%$$
%A(\Te,X)\cong A(\Te',X)\quad \Leftrightarrow\quad  A(\Te)\cong A(\Te') 
%\quad \Leftrightarrow \quad 
%\hat{A}(\Te)\cong \hat{A}(\Te')
%$$
$A(\Te,X)\cong A(\Te',X)$ if and only if $A(\Te)\cong A(\Te')$ (equivalently, 
$\hat{A}(\Te)\cong \hat{A}(\Te')$).

If $A$ has an identity element ${\bf 1}_A$ (i.e. $A=\hat{A}$) then we can and 
will assume that ${\bf 1}_A\in A_i$ for all $i$. Put $B_i=A_i$ if ${\bf 1}_A\in A$
and $B_i=A_i+\bbF {\bf 1}_{\hat{A}}$ otherwise, see the proof of Lemma \ref{unital}.
Then $B_i$ is semisimple, with possibly one extra 1-dimensional simple component.
Moreover, all embeddings $B_i\to B_{i+1}$ are unital and $\hat{A}=\dlim B_i$.
Recall that $X$ is the type of $A$. 
We will denote by   $\mathcal{B}(\Te)$ the Bratteli diagram 
 $\mathcal{B}(\hat{A})$ of the 
algebra $\hat{A}$ with respect to the sequence $B_1\to B_2\to B_3\to\dots$. 

If ${\bf 1}_A\in A$ and the type $X={\bf S},{\bf O}$, then all $z_i=0$ and  it is easy 
to see that $\mathcal{B}(\Te)$ is
\begin{equation}\label{bratS1}
\xymatrix{
{\bullet} \ar@{-}[rr]|<<<<<<<{l_1}  & & 
{\bullet} \ar@{-}[rr]|<<<<<<<{l_2}  & & 
{\bullet} \ar@{-}[rr]|<<<<<<<{l_3}  & & \dots
}
\end{equation}
The locally semisimple algebras of this type are just the limits
of ``pure diagonal'' matrix embeddings $M_{n_i}\to M_{n_{i+1}}$ given by
$M\mapsto \diag(M,\dots,M)$ ($l_i$ blocks), $M\in M_{n_i}(\bbF)$.
They were first classified by Glimm \cite{glimm} (in $\mathbb{C}^*$-algebras setting). 
It is easy to see that two algebras of this type are isomorphic
if and only if their corresponding supernatural numbers 
$\Pi=l_1 l_2 l_3 \dots$ are equal.

If ${\bf 1}_A\not\in A$ and the type $X={\bf S},{\bf O}$, then  $\mathcal{B}(\Te)$ is
\begin{equation}\label{bratS2}
\xymatrix{
{\bullet} \ar@{-}[rr]|<<<<<<<{l_1}  & & 
{\bullet} \ar@{-}[rr]|<<<<<<<{l_2}  & & 
{\bullet} \ar@{-}[rr]|<<<<<<<{l_3}  & & \dots
\\
{\bullet} \ar@{-}[rr]|<<<<<<<{1}  \ar@{-}[urr]|<<<<<<<{z_1} & & 
{\bullet} \ar@{-}[rr]|<<<<<<<{1}  \ar@{-}[urr]|<<<<<<<{z_2} & & 
{\bullet} \ar@{-}[rr]|<<<<<<<{1}  \ar@{-}[urr]|<<<<<<<{z_3} & &  \dots
}
\end{equation}
The corresponding locally semisimple algebras $A(\Te)$ are the direct limits
of matrix embeddings of the shape (\ref{e01}). 
They were first classified by Dixmier \cite{Di} (in $\mathbb{C}^*$-algebras setting). 
Dixmier's parametrization consists of  the 
supernatural number $\Pi=l_1 l_2 l_3 \dots$ and one real parameter 
$\theta$, which is in fact the inverse of our density index $\delta$, see below. 
The diagrams of this shape also parametrize so-called ``diagonal''
direct limits of finite symmetric and alternating groups \cite{lavren1}.

If ${\bf 1}_A\in A$ and the type $X={\bf A}$, then  $\mathcal{B}(\Te)$ is
\begin{equation}\label{bratA1}
\xymatrix{
{\bullet} \ar@{-}[rr]|<<<<<<<{l_1}  \ar@{-}[drr]|<<<<<<<{r_1} & & 
{\bullet} \ar@{-}[rr]|<<<<<<<{l_2}  \ar@{-}[drr]|<<<<<<<{r_2} & & 
{\bullet} \ar@{-}[rr]|<<<<<<<{l_3}  \ar@{-}[drr]|<<<<<<<{r_3} & & \dots
\\
{\bullet} \ar@{-}[rr]|<<<<<<<{l_1}  \ar@{-}[urr]|<<<<<<<{r_1} & & 
{\bullet} \ar@{-}[rr]|<<<<<<<{l_2}  \ar@{-}[urr]|<<<<<<<{r_2} & & 
{\bullet} \ar@{-}[rr]|<<<<<<<{l_3}  \ar@{-}[urr]|<<<<<<<{r_3} & & \dots
}
\end{equation}
The corresponding algebras were first classified by 
Fack and Mar\'echal \cite{fm1}  (in $\mathbb{C}^*$-algebras setting). 

If ${\bf 1}_A\not\in A$ and the type $X={\bf A}$, then  $\mathcal{B}(\Te)$ is
\begin{equation}\label{bratA2}
\xymatrix{
{\bullet} \ar@{-}[rr]|<<<<<<<{l_1}  \ar@{-}[drr]|<<<<<<<{r_1} & & 
{\bullet} \ar@{-}[rr]|<<<<<<<{l_2}  \ar@{-}[drr]|<<<<<<<{r_2} & & 
{\bullet} \ar@{-}[rr]|<<<<<<<{l_3}  \ar@{-}[drr]|<<<<<<<{r_3} & & \dots
\\
{\bullet} \ar@{-}[rr]|<<<<<<<{l_1}  \ar@{-}[urr]|<<<<<<<{r_1} & & 
{\bullet} \ar@{-}[rr]|<<<<<<<{l_2}  \ar@{-}[urr]|<<<<<<<{r_2} & & 
{\bullet} \ar@{-}[rr]|<<<<<<<{l_3}  \ar@{-}[urr]|<<<<<<<{r_3} & & \dots
\\
{\bullet} \ar@{-}[rr]|<<<<<<<{1}  \ar@{-}[uurr]|<<<<<<<{z_1} \ar@{-}[urr]|<<<<<<<{z_1}& &
{\bullet} \ar@{-}[rr]|<<<<<<<{1}  \ar@{-}[uurr]|<<<<<<<{z_2} \ar@{-}[urr]|<<<<<<<{z_2}& &
{\bullet} \ar@{-}[rr]|<<<<<<<{1}  \ar@{-}[uurr]|<<<<<<<{z_3} \ar@{-}[urr]|<<<<<<<{z_3}& & \dots
}
\end{equation}
This is the most general case. We parametrize the corresponding algebras
by two supernatural numbers and two real parameters (see Theorem \ref{main}).
%The diagrams of this shape appear in the classification of ``diagonal''
%direct limits of finite dimensional simple Lie algebras 
% \cite{bzh}. 

Let $B=\dlim B_i$ be a unital locally semisimple
algebra  and let $K_0(B)$ be its  Grothendieck group 
with positive cone $K_0(B)^+$. Note that 
the homomorphism $B_i\to B_{i+1}$ induces the homomorphism of the abelian 
groups $K_0(B_i)\to K_0(B_{i+1})$ and $K_0(B)$ can be obtained as the direct limit
$\dlim K_0(B_i)$. Since $B_i$ are finite dimensional and semisimple, one has
$(K_0(B_i),K_0(B_i)^+)=(\mathbb{Z}^{k_i},\mathbb{Z}_+^{k_i})$ 
where $k_i$ is the number of the simple components of $B_i$. Therefore the abelian group
$K_0(B)$ is the direct limit of the sequence
$$
\mathbb{Z}^{k_1}\to \mathbb{Z}^{k_2} \to \dots \to\mathbb{Z}^{k_i}\to \mathbb{Z}^{k_{i+1}}\to \dots
$$
Moreover, the embedding on the $i$th level is given by the adjacency (or multiplicities) matrix of 
the $i$th level of the Bratteli diagram of $B$. For example, for the algebra $\hat{A}$ in (\ref{bratA2})
the group $K_0(\hat{A})$ is the direct limit of the sequence 
$\mathbb{Z}^3\to \mathbb{Z}^3\to\dots\to \mathbb{Z}^3\to\dots$ 
were the embedding on the $i$th level is given by the matrix 
$
\left(
\begin{array}{ccc}
l_i & r_i & z_i \\
r_i & l_i & z_i \\
0 & 0 & 1
\end{array}
\right)
.$

Let  ${\bf 1}_B$ be  the identity element of $B$ and let $[{\bf 1}_B]$
be the corresponding element of $K_0(B)^+$. 
The triple $(K_0(B),K_0(B)^+,[{\bf 1}_B])$ is called the {\em dimension group} 
of $B$ and is a complete invariant for unital locally semisimple
algebras. More exactly the following is true.

\begin{ttt}\cite{elliot} \label{dimgr}
Let $B_1$ and $B_2$ be unital locally semisimple
algebras. Then $B_1\cong B_2$ if and only if there is
an order-isomorphism $\varphi:K_0(B_1)\to K_0(B_2)$ such that
$\varphi([{\bf 1}_{B_1}])=[{\bf 1}_{B_2}]$.
\end{ttt}

A similar result holds for non-unital algebras if one replaces $[{\bf 1}_B]$ by the {\em scale} of $K_0(B)$.

For a triple sequence $\Te$ we denote by $G(\Te)$ the dimension group
of the unital locally semisimple algebra $\hat{A}(\Te)$.

\section{The classification of algebras of the same type and the corresponding dimension groups} 
\label{same}

In this section, $\Te=(l_i,r_i,z_i)_{i\in \bbN}$ is the triple sequence 
of  the canonically represented locally
involution simple algebra $A=A(\Te,X)$ of type $X$. 
Recall that $l_i\ge r_i$ and $l_i+r_i\ge 1$ for all $i$.
The degrees $n_i$ of the subalgebras $A_i$ satisfy the following: 
$n_1=1$ and $n_{i+1}=(l_i+r_i)n_i+z_i$ for all $i\ge 1$.

Set $s_i=l_i+r_i$, $c_i=l_i-r_i$ ($i=1,2,\dots$),  $\Se=(s_i)_{i\in \bbN}$, 
$\Ce=(c_i)_{i\in \bbN}$, $s_i^k=s_i\dots s_{k-1}$ and  $c_i^k=c_i\dots c_{k-1}$. 
Put $\delta_i=s_1^{i}/n_i$. 
Then
\begin{equation}\label{eqdelta}
\delta_{i+1}=\frac{s_1^{i+1}}{n_{i+1}}=\frac{s_1^{i}s_i}{n_is_i+z_i}=
\frac{s_1^{i}}{n_i+(z_i/s_i)}\le \delta_i.
\end{equation}
The limit 
$$
\delta=\lim_{i\to\infty}\delta_i
$$
is called the {\em density index} of $\Te$ and is denoted by 
$\delta(\Te)$. Since $\delta_2=s_1/n_2=1$, we have $0\le\delta\le1$.
If $\delta=0$, then the triple sequence is called {\em sparse}. If there 
exists
$i$ such that for all $j>i$ we have $\delta_j=\delta_i\ne0$, then
the triple sequence is called {\em pure}. In view of (\ref{eqdelta}) 
this is equivalent
to the following. There exists $i$ such that for all $j\ge i$
we have $z_j=0$. In this case, 
by removing a finite number of terms from the canonically represented sequence
without changing the limit algebra, 
we may and will assume that $z_i=0$ for all $i$. 
We say that the triple sequence is {\em dense} if and only if
$0<\delta<\delta_i$ for all $i$.

If there exists $i$ such that $c_j=s_j$ (equivalently, $r_j=0$)
for all $j\ge i$, then $\Te$ is called {\em one-sided}. Otherwise,
it is called {\em two-sided}. If for each $i$ there exists $j>i$
such that $c_j=0$ (equivalently, $l_j=r_j$), then 
$\Te$ is called (two-sided) {\em symmetric}. Otherwise
it is called {\em non-symmetric}. In the latter case we may and will
assume that $c_i>0$ for all $i\in\bbN$. Set
$\sigma_i=\frac{c_1\dots c_i}{s_1\dots s_i}$. The limit
$$
\sigma=\lim_{i\to\infty}\sigma_i
$$ 
is called the {\em symmetry index} of
$\Te$ and is denoted by $\sigma(\Te)$. Observe that $0\le\sigma\le1$.
Two-sided non-symmetric triple sequences with $\sigma=0$ are called
{\em weakly non-symmetric}, and those with
$\sigma\ne0$ are called {\em strongly non-symmetric}.

Thus all triple sequences can be partitioned into three classes 
with respect to density
and into four classes with respect to symmetry.

\medskip
    {\bf Density types}                      

\begin{itemize}
\item[(D1)]  Sparse ($\delta=0$).                     

\item[(D2)]  Dense ($\delta_i>\delta>0$ for all $i$). 

\item[(D3)]  Pure ($\delta_i=\delta>0$ for some $i$).
\end{itemize}

    \medskip                                                    
    {\bf Symmetry types}                      

\begin{itemize}
\item[(S1)]  One-sided ($r_j=0$ for all $j\gg1$). 

\item[(S2)]  Two-sided symmetric ($l_j=r_j$ for an infinite set of $j$).

\item[(S3)]  Two-sided weakly non-symmetric 
($r_j>0$ for an infinite set of $j$, 
$l_k>r_k$ for all $k\gg1$,
and $\sigma=0$).

\item[(S4)]  Two-sided strongly non-symmetric 
($r_j>0$ for an infinite set of $j$, 
$l_k>r_k$ for all $k\gg1$,
and $\sigma\ne0$).
\end{itemize}

Now we are ready to prove 
our main classification result for algebras of the same type.

\begin{ttt}\label{main}
Let $\Te=\{(l_i,r_i,z_i)\mid {i\in \bbN}\}$ and $\Te'=\{(l_i',r_i',z_i')\mid {i\in \bbN}\}$ be triple sequences 
and let $X={\bf A},{\bf S}$ 
or $\bf O$. 
Set $\delta=\delta(\Te)$, $\sigma=\sigma(\Te)$,
$\delta'=\delta(\Te')$ and $\sigma'=\sigma(\Te')$. Then
the locally involution simple algebras $A(\Te,X)$ and  $A(\Te',X)$
(respectively, the locally semisimple algebras $A(\Te)$ and  $A(\Te')$;
respectively, the dimension groups $G(\Te)$ and  $G(\Te')$)
are isomorphic if and only if the 
following conditions
hold.
\begin{enumerate}
\item[$({\cal A}_1)$] The triple sequences $\Te$ and $\Te'$ have the same 
density 
type.
\item[$({\cal A}_2)$] $\Pi(\Se)\stackrel{\bbQ}{\sim}\Pi(\Se')$.
\item[$({\cal A}_3)$] $\frac{\delta}{\delta'}\in\frac{\Pi(\Se)}{\Pi(\Se')}$
for dense and pure triple sequences (types {\rm(D2)} and {\rm(D3)}).
\item[$({\cal B}_1)$] The triple sequences $\Te$ and $\Te'$ have the same 
symmetry 
type.
\item[$({\cal B}_2)$] $\Pi(\Ce)\stackrel{\bbQ}{\sim}\Pi(\Ce')$ for
two-sided non-symmetric triple sequences (types {\rm(S3)} and {\rm(S4)}).
\item[$({\cal B}_3)$] There exists $\alpha\in\frac{\Pi(\Se)}{\Pi(\Se')}$
such that $\alpha\frac{\sigma}{\sigma'}\in\frac{\Pi(\Ce)}{\Pi(\Ce')}$
for two-sided strongly non-symmetric triple sequences (type {\rm(S4)}).
Moreover, $\alpha=\frac{\delta}{\delta'}$ if in addition the triple 
sequences
are dense or pure (types {\rm(D2)} and {\rm(D3)}).
\end{enumerate}
\end{ttt}

\noindent{\em Proof}. The proof is similar to that in the case of 
Lie algebras in characteristic zero (see \cite{bzh}). First we will prove 
necessity. 
By Theorem \ref{aareduced} and Propositions \ref{adiso} and \ref{dimgr},
it is enough to prove the result for the locally semisimple algebras
$A(\Te)$ and  $A(\Te')$.
We will prove the following more general statement (which will later be 
used for intertype isomorphisms, Theorem \ref{t51}). If 
$A(\Te,X)\cong A(\Te',X')$ (we do not demand that $X=X'$), then $\Te$ and 
$\Te'$ satisfy the conditions $({\cal A}_1),({\cal A}_2),({\cal A}_3)$. 
Moreover, if  
$X=X'=\bf A$, then the conditions $({\cal B}_1),({\cal B}_2),({\cal B}_3)$ hold. 
Let $( A_i)_{i\in I}$ and $( A_j')_{j\in J}$ $(I\cong J\cong \bbN)$ be 
canonically represented sequences of involution simple algebras of types $X$ and $X'$, 
corresponding to the triple sequences $\Te$ and $\Te'$, respectively. 
We have $ A\cong A'$ where 
$ A=\dlim A_i$, $ A'=\dlim A_i'$. By Proposition~\ref{p50}, 
there exist 
subsequences $i_1<i_2<\dots$ of $I$, 
$j_1<j_2<\dots$ of $J$, and embeddings $\ep_k: A_{i_k}\to A_{j_k}'$, 
$\ep_k': A_{j_k}'\to A_{i_{k+1}}$ ($k=1,2,\dots$)  such that the 
following diagram 
is commutative.
\begin{equation}\label{e9}
\begin{array}{ccccccccccccc}
 A_{i_1} & \longrightarrow & \dots & \longrightarrow & 
 A_{i_k} & \longrightarrow &  A_{i_{k+1}} & \longrightarrow & \dots &
\longrightarrow &  A_{i_m} & \longrightarrow & \dots 
\\
{\downarrow}\lefteqn{\scriptstyle \ep_1} &
{\nearrow} \lefteqn{\scriptstyle \ep_1'}& &
{\nearrow} &
{\downarrow}\lefteqn{\scriptstyle \ep_k} &
{\nearrow} \lefteqn{\scriptstyle \ep_k'}& 
{\downarrow}\lefteqn{\scriptstyle \ep_{k+1}} &
{\nearrow} & & {\nearrow} &
{\downarrow}\lefteqn{\scriptstyle \ep_m} &
{\nearrow} & 
\\
 A_{j_1}'& \longrightarrow & \dots & \longrightarrow & 
 A_{j_k}'& \longrightarrow &  A_{j_{k+1}}'& \longrightarrow & \dots &
\longrightarrow &  A_{j_m}'& \longrightarrow & \dots 
\end{array}
\end{equation}
Let 
$(p_k,q_k,u_k)$ 
(resp., $(p_k',q_k',u_k')$) be the signature of $\ep_k$ (resp., $\ep_k'$). 
Let $n_i$ be the degree of $A_i$. Set $s_i=l_i+r_i$, $c_i=l_i-r_i$, $\delta_i=s_1^i/n_i$, 
$\de=\lim_{i\to\infty}\de_i$. The numbers $n_j'$, $s_j'$,... for the algebra $A'$ are defined similarly. We have
\begin{equation}\label{e10}
n_{j_m}'=(p_m+q_m)n_{i_m}+u_m=(p_m+q_m)s_1^{i_m}\de_{i_m}^{-1}+u_m=
(p_m+q_m)s_1^{i_k}s_{i_k}^{i_m}\de_{i_m}^{-1}+u_m.
\end{equation}
On the other hand,
\begin{equation}\label{e11}
n_{j_m}'=s_1'^{j_m}(\de_{j_m}')^{-1}=
s_1'^{j_k}s_{j_k}'^{j_m}(\de_{j_m}')^{-1}.
\end{equation}
In view of commutativity of the diagram and by Corollary~\ref{c80} we have
\begin{equation}\label{e11a}
s_{i_k}^{i_m}(p_m+q_m)=(p_k+q_k)s_{j_k}'^{j_m}.
\end{equation}
Dividing (\ref{e10}) and (\ref{e11}) by $s_{j_k}'^{j_m}$, we get
$
(p_k+q_k)s_1^{i_k}\de_{i_m}^{-1}+u_m/s_{j_k}'^{j_m}=
s_1'^{j_k}(\de_{j_m}')^{-1},
$
so 
\begin{equation}\label{e12}
(p_k+q_k)s_1^{i_k}\de_{j_m}'\le s_1'^{j_k}\de_{i_m}.
\end{equation}
Taking $m\to\infty$, we obtain $(p_k+q_k)s_1^{i_k}\de'\le s_1'^{j_k}\de$. 
Similarly, we get $(p_k'+q_k')s_1'^{j_k}\de\le s_1^{i_{k+1}}\de'$. By 
Corollary~\ref{c80}, we have $(p_k+q_k)(p_k'+q_k')=s_{i_k}^{i_{k+1}}$. Hence
$$
(p_k+q_k)s_1^{i_k}\de'\le s_1'^{j_k}\de\le 
(p_k'+q_k')^{-1}s_1^{i_{k+1}}\de'=(p_k+q_k)s_1^{i_k}\de'.
$$
Therefore
\begin{eqnarray}
(p_k+q_k)s_1^{i_k}\de'&=&s_1'^{j_k}\de,       \label{e13}\\
(p_k'+q_k')s_1'^{j_k}\de&=&s_1^{i_{k+1}}\de'. \label{e14}
\end{eqnarray}
Clearly $\de=0$ if and only if $\de'=0$. Therefore $\Te$ is sparse if and 
only if $\Te'$ is so. If the triple sequence $\Te$ is pure, then 
$\de=\de_{i_m}$ for some $m$. 
Subtracting (\ref{e13}) from (\ref{e12}), we get
$$
0\le (p_k+q_k)s_1^{i_k}(\de_{j_m}'-\de')\le s_1'^{j_k}(\de_{i_m}-\de)=0.
$$
Therefore $\de_{j_m}'=\de'$, so $\Te'$ is also pure. By symmetry, $\Te$ is 
pure if and only if $\Te'$ is pure. So $({\cal A}_1)$ holds.

By (\ref{e11a}), $s_{i_k}^{i_m}$ divides $(p_k+q_k)s_{j_k}'^{j_m}$ for all 
$m>k$. On the other hand, in view of commutativity of the diagram we have
\begin{equation}\label{e15}
s_{i_k}^{i_{m+1}}=(p_k+q_k)s_{j_k}'^{j_m}(p_m'+q_m'),
\end{equation}
so $(p_k+q_k)s_{j_k}'^{j_m}$ divides $s_{i_k}^{i_{m+1}}$. Therefore by 
Proposition~\ref{p33},
\begin{equation}\label{e15a}
\Pi(\Se_{i_k})=(p_k+q_k)\Pi(\Se_{j_k}'),
\end{equation}
where $\Se_{i_k}=(s_{i_k},s_{i_k+1},\dots)$, 
$\Se_{j_k}'=(s_{j_k}',s_{j_k+1}',\dots)$. It follows that 
$\Pi(\Se)\stackrel{\bbQ}{\sim}\Pi(\Se')$, so $({\cal A}_2)$ holds.

Finally, if $\de$ and $\de'$ are nonzero (dense or pure sequences), then by 
(\ref{e13}) and (\ref{e14}), $s_1^{i_k}$ divides $(\de/\de')s_1'^{j_k}$ and 
$(\de/\de')s_1'^{j_k}$ divides $s_1^{i_{k+1}}$ for any $k$. Therefore by 
Proposition~\ref{p33}, $\Pi(\Se)=(\de/\de')\Pi(\Se')$, and $({\cal A}_3)$ 
holds.

Assume now that $X=X'=\bf A$. By Corollary~\ref{c80}, one can write down 
equalities for ``differences'' similar to (\ref{e11a}) and (\ref{e15}):
\begin{equation}\label{e16}
c_{i_k}^{i_m}(p_m-q_m)=(p_k-q_k)c_{j_k}'^{j_m};
\end{equation}
\begin{equation}\label{e17}
c_{i_k}^{i_{m+1}}=(p_k-q_k)c_{j_k}'^{j_m}(p_m'-q_m').
\end{equation}
If $\Te'$ is symmetric, then by definition, for each $k$ there exists $m$ such 
that $c_{j_k}'^{j_m}=0$. It follows from (\ref{e17}) that 
$c_{i_k}^{i_{m+1}}=0$, so $\Te$ is symmetric. Therefore, $\Te$ is symmetric 
if and only if $\Te'$ is so. Assume that $\Te$ is non-symmetric. Recall 
that in this case one can suppose that all $c_i$ and $c_j'$ are nonzero. 
Dividing (\ref{e17}) by (\ref{e15}), we get 
\begin{equation}\label{e18}
\frac{c_{i_k}^{i_{m+1}}}{s_{i_k}^{i_{m+1}}}=
\frac{(p_k-q_k)}{(p_k+q_k)}
\frac{c_{j_k}'^{j_m}}{s_{j_k}'^{j_m}}
\frac{(p_m'-q_m')}{(p_m'+q_m')},
\end{equation}
or equivalently,
\begin{equation}\label{e18a}
\sigma_1^{i_{m+1}}\cdot
\frac{s_1^{i_k}}{c_1^{i_k}}=
\sigma_1'^{j_m}\cdot
\frac{(p_k-q_k)}{(p_k+q_k)}
\frac{s_1'^{j_k}}{c_1'^{j_k}}
\frac{(p_m'-q_m')}{(p_m'+q_m')}.
\end{equation}
Taking $m\to\infty$, we get 
\begin{equation}\label{e19}
\sigma\cdot
\frac{s_1^{i_k}}{c_1^{i_k}}\le
\sigma'\cdot
\frac{(p_k-q_k)}{(p_k+q_k)}
\frac{s_1'^{j_k}}{c_1'^{j_k}}.
\end{equation}
Similarly, dividing (\ref{e16}) by (\ref{e11a}) and taking $m\to\infty$, we 
get
\begin{equation}\label{e20}
\sigma\cdot
\frac{s_1^{i_k}}{c_1^{i_k}}\ge
\sigma'\cdot
\frac{(p_k-q_k)}{(p_k+q_k)}
\frac{s_1'^{j_k}}{c_1'^{j_k}}.
\end{equation}
Combining with (\ref{e19}), we obtain
\begin{equation}\label{e21}
\sigma\cdot
\frac{s_1^{i_k}}{c_1^{i_k}}=
\sigma'\cdot
\frac{(p_k-q_k)}{(p_k+q_k)}
\frac{s_1'^{j_k}}{c_1'^{j_k}}.
\end{equation}
It follows that $\sigma=0$ if and only if $\sigma'=0$. That is, $\Te$ is 
weakly non-symmetric if and only if $\Te'$ is so. Assume that $\Te'$ is 
one-sided. Then $\sigma'=\sigma_1'^{j_m}$ for some $m$. Subtracting 
(\ref{e21}) from (\ref{e18a}), we have 
$0\le(\sigma_1^{i_{m+1}}-\sigma)s_1^{i_k}/c_1^{i_k}\le0$. Therefore 
$\sigma_1^{i_{m+1}}=\sigma$, i.e. $\Te$ is one-sided. So $({\cal B}_1)$
holds.
 
Similarly to (\ref{e15a}), one can get
\begin{equation}\label{e21a}
\Pi(\Ce_{i_k})=(p_k-q_k)\Pi(\Ce_{j_k}').
\end{equation}
It follows that 
$\Pi(\Ce)\stackrel{\bbQ}{\sim}\Pi(\Ce')$, so $({\cal B}_2)$ holds.

Assume now that $\Te$ and $\Te'$ are strongly non-symmetric, i.e.  
$\sigma\ne 0$ and $\sigma'\ne0$. Set $\alpha=(p_k+q_k)s_1^{i_k}/s_1'^{j_k}$. 
Then (\ref{e21}) can be rewritten in the form
\begin{equation}\label{e22}
\frac{\sigma}{\sigma'}\alpha c_1'^{j_k}=(p_k-q_k)c_1^{i_k}.
\end{equation}
Observe that $\alpha\in\frac{\Pi(\Se)}{\Pi(\Se')}$. Indeed, using 
(\ref{e15a}), we have
$$
\alpha\Pi(\Se')=(p_k+q_k)s_1^{i_k}\Pi(\Se_{j_k}')=
s_1^{i_k}\Pi(\Se_{i_k})=\Pi(\Se).
$$
Moreover, if $\Te$ and $\Te'$ are dense or pure, then by (\ref{e13}), 
$\alpha=\delta/\delta'$. It follows from (\ref{e22}) and (\ref{e21a}) that
$$
\frac{\sigma}{\sigma'}\alpha\Pi(\Ce')=
(p_k-q_k)c_1^{i_k}\Pi(\Ce_{j_k}')=
c_1^{i_k}\Pi(\Ce_{i_k})=\Pi(\Ce).
$$
Therefore, 
$\frac{\sigma}{\sigma'}\alpha\in\frac{\Pi(\Ce)}{\Pi(\Ce')}$.
This proves $({\cal B}_3)$. \hfill$\Box$\par\medskip

To prove the sufficiency in Theorem~\ref{main}, we need the following 
lemma.

\begin{lll}\label{l17}
Let $\Te$ and $\Te'$ satisfy the conditions 
$({\cal A}_1),({\cal A}_2),({\cal A}_3)$, 
$({\cal B}_1),({\cal B}_2)$, $({\cal B}_3)$ of the theorem. Fix 
$\alpha\in\frac{\Pi(\Se)}{\Pi(\Se')}$ ($\alpha=\delta/\delta'$ if $\Te$ 
and $\Te'$ are dense or pure), $\beta\in\frac{\Pi(\Ce)}{\Pi(\Ce')}$ for 
the case of two-sided non-symmetric triple sequences 
($\beta/\alpha=\sigma/\sigma'$ if $\Te$ and $\Te'$ are strongly 
non-symmetric). Let $i,j,a,b$ be integers such that

$(a)\quad\alpha s_1'^j=a s_1^i$, 

$(b)\quad\beta c_1'^j=b c_1^i$ 
(for two-sided non-symmetric $\Te$ and $\Te'$).

\noindent Then there exists $k>i$ such that $a'=s_i^k/a$ and $b'=c_i^k/b$ are 
integers of the same parity 
($a'$ is even and $c_i^k=0$ for the case of symmetric $\Te$ and $\Te'$), 
$a'\ge b'$ and $n_k\ge a' n_j'$.
\end{lll}

\begin{pf}
If otherwise is not specified we assume 
that $\Te$ and $\Te'$ are two-sided non-symmetric. The case 
of one-sided and symmetric sequences can be settled by removing from the 
proof the arguments with $c_i$, $\beta$, $b$. 

Since $\alpha\in\frac{\Pi(\Se)}{\Pi(\Se')}$ and $\alpha s_1'^j=a s_1^i$, 
we have
\begin{equation}\label{e30}
\Pi(\Se_i)=\alpha(s_1^i)^{-1}s_1'^j\Pi(\Se_j')=a\Pi(\Se_j').
\end{equation}
Similarly, we get 
\begin{equation}\label{e31}
\Pi(\Ce_i)=b\Pi(\Ce_j').
\end{equation}
Therefore there exists $k_1>i$ such that $a'=s_i^k/a$ and $b'=c_i^k/b$ are 
integers for all $k\ge k_1$. Since for each $m$ the integers 
$s_m'=l_m'+r_m'$ and $c_m'=l_m'-r_m'$ have the same parity, 2 divides 
$\Pi(\Se_j')$ if and only if 2 divides $\Pi(\Ce_j')$ (for symmetric 
sequences 2 divides $\Pi(\Se_j')$ always). Therefore by (\ref{e30}) and 
(\ref{e31}), there exists $k_2\ge k_1$ such that the integers $a'$ and $b'$ 
have the same parity 
($a'$ is even and $c_i^k=0$ for the case of symmetric $\Te$ and $\Te'$) 
for all $k\ge k_2$. Set $\gamma_k=b'/a'$. In view of $(a)$ and $(b)$, 
we have 
$$
\gamma_k=\frac{c_i^k}{b}\cdot\frac{a}{s_i^k}=
\frac{c_1^ic_i^k}{\beta c_1'^j}\cdot\frac{\alpha s_1'^j}{s_1^is_i^k}= 
\frac{\alpha}{\beta}\cdot\frac{\sigma_1^k}{\sigma_1'^{j}}.
$$
If $\Te$ and $\Te'$ are weakly non-symmetric, then $\sigma_1^k\to0$ as 
$k\to\infty$, so $\gamma_k\to0$. 
If $\Te$ and $\Te'$ are strongly non-symmetric, then 
by assumption $\beta/\alpha=\sigma/\sigma'$, so
$$
\gamma_k\to\frac{\alpha}{\beta}\cdot\frac{\sigma}{\sigma_1'^j}=
\frac{\sigma'}{\sigma_1'^j}<1
$$
as $k\to\infty$. 
In both cases there exists $k_3\ge k_2$ such that $\gamma_k\le 1$ (i.e. 
$a'\ge b'$) for all $k\ge k_3$.

Set $\nu_k=n_k/a'-n_j'$. We have to show that $\nu_k\ge0$
for sufficiently large $k$. One has
$$
\nu_k=\frac{n_k}{a'}-n_j'=
\frac{s_1^k}{a'\delta_k}-\frac{s_1'^j}{\delta_j'}=
\frac{a s_1^i}{\delta_k}-\frac{s_1'^j}{\delta_j'}=
s_1'^j(\frac{\alpha}{\delta_k}-\frac{1}{\delta_j'}).
$$
(The last equality follows from $(a)$.)
If $\Te$ and $\Te'$ are sparse, then $\delta_k\to0$ as $k\to\infty$, so 
$\nu_k\to+\infty$. Therefore there exists $k_4\ge k_3$ such that $\nu_k\ge0$
for all $k\ge k_4$. 
Let $\Te$ and $\Te'$ be dense. Then $\alpha=\delta/\delta'$ and 
$\delta_j'>\delta'$. Therefore 
$$
\nu_k=s_1'^j(\frac{\delta}{\delta_k}\cdot
\frac{1}{\delta'}-\frac{1}{\delta_j'})\to
s_1'^j(\frac{1}{\delta'}-\frac{1}{\delta_j'})>0,
$$
as $k\to\infty$. Hence there exists $k_4\ge k_3$ such that $\nu_k\ge0$
for all $k\ge k_4$.  
Let $\Te$ and $\Te'$ be pure. Then there exists $k_4\ge k_3$ such that 
$\delta=\delta_k$ for all $k\ge k_4$. Therefore
$$
\nu_k=s_1'^j(\frac{1}{\delta'}-\frac{1}{\delta_j'})\ge0,
$$
for all $k\ge k_4$. So each $k\ge k_4$ satisfies the assumptions of the theorem.
\end{pf}

{\em Proof of sufficiency in Theorem~\ref{main}}. According to 
Proposition~\ref{p50} we have to construct 
sequences $i_1<i_2<\dots$, 
$j_1<j_2<\dots$, and embeddings $\ep_k: A_{i_k}\to A_{j_k}'$, 
$\ep_k': A_{j_k}'\to A_{i_{k+1}}$ ($k=1,2,\dots$)  such that the 
diagram~(\ref{e9})  
is commutative. 
Fix 
$\alpha\in\frac{\Pi(\Se)}{\Pi(\Se')}$ ($\alpha=\delta/\delta'$ if $\Te$ 
and $\Te'$ are dense or pure) and $\beta\in\frac{\Pi(\Ce)}{\Pi(\Ce')}$ for 
the case of two-sided non-symmetric triple sequences 
($\beta/\alpha=\sigma/\sigma'$ if $\Te$ and $\Te'$ are strongly 
non-symmetric). Fix also $j_0\in J$. Since 
$\Pi(\Se')=\alpha^{-1}\Pi(\Se)$ and $\Pi(\Ce')=\beta^{-1}\Pi(\Ce)$, by 
Proposition~\ref{p33}, there exists $i_1\in I$ such that

$(a_0)\quad\alpha^{-1} s_1^{i_1}=a_0 s_1'^{j_0}$,

$(b_0)\quad\beta^{-1} c_1^{i_1}=b_0 c_1'^{j_0}$ (for two-sided 
non-symmetric $\Te$ and $\Te'$)

\noindent 
where $a_0,b_0\in\bbN$. Applying Lemma~\ref{l17} (interchanging $\Te$ and 
$\Te'$), we find $j_1$ such that $a_1=s_{j_0}'^{j_1}/a_0$ and 
$b_1=c_{j_0}'^{j_1}/b_0$ are integers of the same parity ($a_1$ is even if 
$\Te$ and $\Te'$ are symmetric), $a_1\ge b_1$ and $n_{j_1}'\ge a_1n_{i_1}$. 
Set $p_1=(a_1+b_1)/2$, $q_1=(a_1-b_1)/2$, $u_1=n_{j_1}'-a_1n_{i_1}$ 
($p_1=q_1=a_1/2$ for symmetric sequences). Consider the canonical embedding 
$\ep_1: A_{i_1}\to A_{j_1}'$ with the signature $(p_1,q_1,u_1)$. We 
have 

$(a_1)\quad\alpha s_1'^{j_1}=a_1 s_1^{i_1}$, 

$(b_1)\quad\beta c_1'^{j_1}=b_1 c_1^{i_1}$.

\noindent Proceed by induction. Assume that 
sequences $i_1<\dots<i_k$, 
$j_1<\dots<j_k$ and embeddings $\ep_1,\ep_1',\dots,\ep_k$ have been 
constructed, and the following conditions hold.

$(a_k)\quad\alpha s_1'^{j_k}=a_k s_1^{i_k}$,

$(b_k)\quad\beta c_1'^{j_k}=b_k c_1^{i_k}$

\noindent 
where $a_k=p_k+q_k$, $b_k=p_k-q_k$. Construct an embedding $\ep_k'$ as follows. By 
Lemma~\ref{l17}, there exists $i_{k+1}>i_k$ such that 
$a_k'=s_{i_k}^{i_{k+1}}/a_k$ and $b_k'=c_{i_k}^{i_{k+1}}/b_k$ are integers 
of the same parity ($a_k'$ is even if $\Te$ and $\Te'$ are symmetric), 
$a_k'\ge b_k'$ and $n_{i_{k+1}}\ge a_k'n_{j_k}'$. Set $p_k'=(a_k'+b_k')/2$, 
$q_k'=(a_k'-b_k')/2$, $u_k'=n_{i_{k+1}}-a_k'n_{j_k}'$ ($p_k'=q_k'=a_k'/2$ 
for symmetric sequences). Since
\begin{eqnarray*}
(p_k+q_k)(p_k'+q_k')&=&a_ka_k'=s_{i_k}^{i_{k+1}}, \\
(p_k-q_k)(p_k'-q_k')&=&b_kb_k'=c_{i_k}^{i_{k+1}}, 
\end{eqnarray*}
and $u_k'\ge0$, by Lemma~\ref{l25} there exists an embedding 
$\ep_k': A_{j_k}'\to A_{i_{k+1}}$ such that $\iota_k=\ep_k'\ep_k$ where 
$\iota_k$ denotes the embedding $ A_{i_k}\to A_{i_{k+1}}$. Observe 
that

$(a_k')\quad\alpha^{-1} s_1^{i_{k+1}}=a_k' s_1'^{j_k}$,

$(b_k')\quad\beta^{-1} c_1^{i_{k+1}}=b_k' s_1'^{j_k}$.

\noindent
Therefore Lemma~\ref{l17} can be applied once more (interchanging $\Te$ and 
$\Te'$). So the result follows by induction.\hfill$\Box$\par\medskip

\begin{rrr}
It is not difficult to see that for pure triple sequences one can always 
assume that all $z_i=0$ (by removing a finite number of terms in the 
sequences). In this case $\delta=1$, so the condition $({\cal A}_3)$ can be rewritten in the form 
$\Pi(\Se)=\Pi(\Se')$.
\end{rrr}

\section{Isomorphisms of algebras of different types}
\label{intertype}

In this section we find conditions under which $A(\Te,X)\cong A(\Te',X')$ 
where $\Te$ and $\Te'$ are triple sequences and $X\ne X'$. We also give a 
general parametrization of countable locally involution simple 
algebras.

\begin{lll}\label{l52}
Let $\Te$ be a two-sided symmetric triple sequence, $\Se=\Se(\Te)$. Then 
$2^{\infty}$ divides $\Pi(\Se)$.
\end{lll}

\begin{pf}
By definition, $l_i=r_i$ (in particular, $s_i=l_i+r_i$ is even) for an 
infinite set of $i$. Therefore $2^{\infty}$ divides $\Pi(\Se)$.
\end{pf}

\begin{ttt}\label{t51}
Let $\Te$, $\Te'$ be triple sequences.
\begin{enumerate}
\item[(i)] Let $\chf\ne2$. Then 
$A(\Te,{\bf A})\cong A(\Te',{\bf O})$ (resp.,$A(\Te,{\bf A})\cong A(\Te',{\bf S})$) if and 
only if $\Te$ is two-sided symmetric, $2^{\infty}$ divides 
$\Pi(\Se')$ and the conditions 
$({\cal A}_1),({\cal A}_2),({\cal A}_3)$ of Theorem~\ref{main} hold.
\item[(ii)] Let $\chf\ne2$.  $A(\Te,{\bf O})\cong A(\Te',{\bf S})$ if and only if 
$2^{\infty}$ divides both $\Pi(\Se)$ and $\Pi(\Se')$,
and the conditions 
$({\cal A}_1),({\cal A}_2),({\cal A}_3)$ of Theorem~\ref{main} hold.
\item[(iii)]  Let $\chf=2$. If $X,X'\in\{{\bf A},{\bf O},{\bf S}\}$
are different, then $A(\Te,X)$ is not isomorphic to $A(\Te',X')$.
\end{enumerate}
\end{ttt}

\begin{pf}
$(i)$. Set $ A=A(\Te,{\bf A})$ and $ A'=A(\Te',{\bf O})$. 
Assume that $ A\cong A'$. Then, as it was established in the proof of Theorem~\ref{main}, 
the conditions $({\cal A}_1)$, $({\cal A}_2)$ and $({\cal A}_3)$ hold. 
Now denote by 
$(x_k,y_k,z_k)$ the signature of $ A_{i_k}\to A_{i_{k+1}}$ (see 
diagram~(\ref{e9})). Since the diagram is commutative, we have 
$x_k=p_kp_k'$ and $y_k=q_kp_k'$ where $(p_k,q_k,u_k)$ and $(p_k',0,u_k')$ 
are the signatures of $\ep_k$ and $\ep_k'$, respectively. By 
Proposition~\ref{p3}$(i)$, $p_k=q_k$, so $x_k=y_k$. Therefore 
$c_{i_k}^{i_{k+1}}=x_k-y_k=0$, so $\Te$ is two-sided symmetric.
By Lemma~\ref{l52}, $2^{\infty}$ divides $\Pi(\Se)$.
Therefore in view of condition $({\cal A}_2)$, 
$2^{\infty}$ divides $\Pi(\Se')$.

{\em Conversely}. Let $ A=A(\Te,{\bf A})$ and  
$ A'=A(\Te',{\bf O})$ be such that $\Te$ is two-sided symmetric,
$2^{\infty}$ divides $\Pi(\Se')$ and the conditions 
$({\cal A}_1)$, $({\cal A}_2)$ and $({\cal A}_3)$ 
hold.  
Then there exists a sequence of indices $j_1<j_2<\dots$ such that 
$s_{j_k}'^{j_{k+1}}$ is even for all $k=1,2,\dots$. By 
Proposition~\ref{p3}$(iii)$, there exists an algebra $ A_k''$ of type $\bf A$ 
and representable embeddings $ A_{j_k}'\to A_k''$ and 
$ A_k''\to A_{j_{k+1}}'$ such that the diagram
$$
\begin{array}{ccccc}
 A_{j_k}'&&\longrightarrow&& A_{j_{k+1}}'\\
           &\searrow&&\nearrow&          \\
                && A_k''&&
\end{array}
$$
is commutative. Set $ A''=\dlim A_k''$. Let $\Te''$ be the 
corresponding triple sequence. We have $ A''=A(\Te'',{\bf A})$. By construction, 
$ A''\cong A'$. Moreover, by the above arguments (the proof of necessity) 
$\Te''$ is symmetric and the conditions 
$({\cal A}_1)$, $({\cal A}_2)$ and $({\cal A}_3)$ 
(for $\Te'$ and $\Te''$) hold. Since the same is true for the pair $\Te$, 
$\Te'$, we conclude that the pair $\Te$, $\Te''$ also satisfies these 
conditions. Indeed, $({\cal A}_1)$ trivially holds. Further, since
$\Pi(\Se')\stackrel{\bbQ}{\sim}\Pi(\Se'')$ and 
$\Pi(\Se)\stackrel{\bbQ}{\sim}\Pi(\Se')$, we have
$\Pi(\Se)\stackrel{\bbQ}{\sim}\Pi(\Se'')$. Finally, if 
$\frac{\delta'}{\delta''}\in\frac{\Pi(\Se')}{\Pi(\Se'')}$ and 
$\frac{\delta}{\delta'}\in\frac{\Pi(\Se)}{\Pi(\Se')}$, then
$$
\Pi(\Se')=\frac{\delta'}{\delta''}\Pi(\Se'')=
\frac{\delta'}{\delta}\Pi(\Se),
$$
so $\frac{\delta}{\delta''}\in\frac{\Pi(\Se)}{\Pi(\Se'')}$. 
Consequently, by Theorem~\ref{main}, $A(\Te,{\bf A})\cong A(\Te'',{\bf A})$, i.e. 
$ A\cong A''$. Therefore $ A\cong A'$. The proof for the case 
$ A'=A(\Te',{\bf S})$ is similar.

$(ii)$. Let $A(\Te,{\bf O})\cong A(\Te',{\bf S})$. Using 
Proposition~\ref{p3}~$(ii),(iii)$, it is not difficult to construct an 
algebra $A(\Te'',{\bf A})\cong A(\Te,{\bf O})\cong A(\Te',{\bf S})$. The claim now follows from 
Theorem~\ref{t51}~$(i)$. To prove the converse statement we construct 
$A(\Te'',{\bf A})$ isomorphic to $A(\Te,{\bf O})$ and use Theorem~\ref{t51}~$(i)$. 

$(iii)$. By definition of 
$A(\Te,X)$, all the corresponding embeddings are representable.
Thus the claim follows from Proposition~\ref{ch2}.
\end{pf}

It remains to discuss the general parametrization. Let $ A$ be 
a locally involution simple associative algebra over $\bbF$ 
of countable dimension. Then by Theorem \ref{rt},
$A$ is canonically representable, i.e. 
$A$ is the direct limit of a 
 sequence $( A_i)_{i\in\bbN}$ of subalgebras of the 
same type $X=\bf A$, $\bf O$, or $\bf S$ such that all embeddings
are canonical. Fix any such system of subalgebras. This gives the
triple sequence $\Te=((l_i,r_i,z_i))_{i\in\bbN}$, the 
sequences of ``sums'' $\Se=(l_i+r_i)_{i\in\bbN}$ and 
(for $X=\bf A$ only) ``differences'' 
$\Ce=(l_i-r_i)_{i\in\bbN}$. Now we can determine the 
density type D$=$(D1), (D2) or (D3), the density index 
$\delta=\delta(\Te)$, supernatural number $\Pi_{\Se}=\Pi(\Se)$, and (for $X=\bf A$ 
only) the symmetry type S$=$(S1), (S2), (S3), or (S4), the symmetry index 
$\sigma=\sigma(\Te)$, and supernatural number $\Pi_{\Ce}=\Pi(\Ce)$. So one can 
associate with any algebra $ A$ a tuple 
$$
{\cal P}( A)=(X,{\rm D,S},\delta,\sigma,\Pi_{\Se},\Pi_{\Ce})
$$
where 
$X,{\rm D,S}$ describe a type of $ A$; $\delta$ and $\sigma$ are real 
numbers ($0\le\delta,\sigma\le1$); $\Pi_{\Se}$ and $\Pi_{\Ce}$ are supernatural 
numbers. For $X={\bf S},{\bf O}$ (and $X=\bf A$ with one-sided or symmetric $\Te$) we use a 
shorter list of invariants: 
$$
 A\mapsto(X,{\rm D},\delta,\Pi_{\Se}).
$$
By Theorem~\ref{main}, the tuples associated with two nonisomorphic algebras are 
distinct. The question under what conditions $ A$ and $ A'$ with 
tuples ${\cal P}( A)$ and ${\cal P}( A')$ are isomorphic has been
resolved in
Theorems~\ref{main} and \ref{t51}.

\end{document}